\def\spin{{\rm Spin}}
\def\endo{{\rm End}}
\def\mult{{\rm mult}}
\def\nmp{{\rm Nm}_p}

\def\cctil{{\widetilde C}}
\def\ftil{{\widetilde F}}
\def\pitil{{\widetilde \pi}}
\def\dual{^{\vee}}
\def\ctil{\cctil_{\eta}}

\def\ahat{{\widehat A}_2}

\def\xieta{{\Xi_{\eta}}}
\def\peta{P_{\eta}}
\def\jtil{{\widetilde J^{2g-2}}}
\def\nm{{\rm Nm}}
\def\pic{{\rm Pic}}

\def\th{K^{1\over 2}}
\def\som{{\bf so}_m}
\def\gab{g_{\alpha \beta}}
\def\xab{\xi_{\alpha \beta}}
\def\phab{\phi_{\alpha \beta}}
\def\eab{\eta_{\alpha \beta}}

\def\ee{{\cal E}}
\def\pp{{\cal P}}

\def\nn{{\cal N}}
\def\mm{{\cal M}}

\def\ll{{\cal L}}

\def\ss{{\cal S}}
\def\ll{{\cal L}}

\def\dd{{\cal D}}
\def\oo{{\cal O}}
\def\uu{{\cal U}}
\def\ww{{\cal W}}

\def\su{{\cal SU}_C}

\def\ses#1#2#3{0\rightarrow {#1}\rightarrow {#2}\rightarrow {#3}\rightarrow 0}

\def\z{{\bf Z}}
\def\n{{\bf N}}
\def\q{{\bf Q}}
\def\c{{\bf C}}
\def\bp{{\bf P}}

\def\ad{{\rm ad}\ }

\def\spec{{\rm Spec \ }}
\def\im{{\rm image \ }}

\def\ker{{\rm ker \ }}

\def\map#1{\ \smash{\mathop{\longrightarrow}\limits^{#1}}\ }

\def\pf{{\it Proof.}\ }

\def\<{\langle}
\def\>{\rangle}

\newcommand{\qed}{{\unskip\nobreak\hfill\hbox{ $\Box$}\medskip\par}}

\font\dynkfont=cmsy10 scaled\magstep4    \skewchar\dynkfont='60
\def\dynk{\textfont2=\dynkfont}
\def\hr#1,#2;{\dimen0=.4pt\advance\dimen0by-#2pt
              \vrule width#1pt height#2pt depth\dimen0}
\def\vr#1,#2;{\vrule height#1pt depth#2pt}
\def\blb#1#2#3#4#5
            {\hbox{\ifnum#2=0\hskip11.5pt
                   \else\ifnum#2=1\hr13,5;\hskip-1.5pt
                   \else\ifnum#2=2\hr13.5,6.5;\hskip-13.5pt
                                  \hr13.5,3.5;\hskip-2pt
                   \else\ifnum#2=3\hr13.7,8;\hskip-13.7pt
                                  \hr13,5;\hskip-13pt
                                  \hr13.7,2;\hskip-2.2pt\fi\fi\fi\fi
                   $#1$
                   \ifnum#4=0
                   \else\ifnum#4=1\hskip-9.2pt\vr22,-9;\hskip8.8pt
                   \else\ifnum#4=2\hskip-10.9pt\vr22,-8.75;\hskip3pt
                                  \vr22,-8.75;\hskip7.1pt
                   \else\ifnum#4=3\hskip-12.6pt\vr22,-8.5;\hskip3pt
                                  \vr22,-9;\hskip3pt
                                  \vr22,-8.5;\hskip5.4pt\fi\fi\fi\fi
                   \ifnum#5=0
                   \else\ifnum#5=1\hskip-9.2pt\vr1,12;\hskip8.8pt
                   \else\ifnum#5=2\hskip-10.9pt\vr1.25,12;\hskip3pt
                                  \vr1.25,12;\hskip7.1pt
                   \else\ifnum#5=3\hskip-12.6pt\vr1.5,12;\hskip3pt
                                  \vr1,12;\hskip3pt
                                  \vr1.5,12;\hskip5.4pt\fi\fi\fi\fi
                   \ifnum#3=0\hskip8pt
                   \else\ifnum#3=1\hskip-5pt\hr13,5;
                   \else\ifnum#3=2\hskip-5.5pt\hr13.5,6.5;
                                  \hskip-13.5pt\hr13.5,3.5;
                   \else\ifnum#3=3\hskip-5.7pt\hr13.7,8;
                                  \hskip-13pt\hr13,5;
                                  \hskip-13.7pt\hr13.7,2;\fi\fi\fi\fi

                   }}
\def\blob#1#2#3#4#5#6#7{\hbox
{$\displaystyle\mathop{\blb#1#2#3#4#5 }_{#6}\sp{#7}$}}
\def\up#1#2{\dimen1=33pt\multiply\dimen1by#1\hbox{\raise\dimen1\rlap{#2}}}
\def\uph#1#2{\dimen1=17.5pt\multiply\dimen1by#1\hbox{\raise\dimen1\rlap{#2}}}
\def\dn#1#2{\dimen1=33pt\multiply\dimen1by#1\hbox{\lower\dimen1\rlap{#2}}}
\def\dnh#1#2{\dimen1=17.5pt\multiply\dimen1by#1\hbox{\lower\dimen1\rlap{#2}}}

\def\rlbl#1{\kern-8pt\raise3pt\hbox{$\scriptstyle #1$}}
\def\llbl#1{\raise3pt\llap{\hbox{$\scriptstyle #1$\kern-8pt}}}
\def\elbl#1{\kern3pt\lower4.5pt\hbox{$\scriptstyle #1$}}
\def\lelbl#1{\rlap{\hbox{\kern-9pt\raise2.5pt\hbox{{$\scriptstyle #1$}}}}}

\def\whtd#1#2#3#4#5{\blob\circ#1#2#3#4{#5}{}}
\def\blkd#1#2#3#4#5{\blob\bullet#1#2#3#4{#5}{}}
\def\whtu#1#2#3#4#5{\blob\circ#1#2#3#4{}{#5}}
\def\blku#1#2#3#4#5{\blob\bullet#1#2#3#4{}{#5}}
\def\whtr#1#2#3#4#5{\blob\circ#1#2#3#4{}{}\rlbl{#5}}
\def\blkr#1#2#3#4#5{\blob\bullet#1#2#3#4{}{}\rlbl{#5}}

\def\rwng{\hbox{$\vbox{\offinterlineskip{
  \hbox{\phantom{}\kern6pt{$\circ$}}\kern-2.5pt\hbox{$\Biggr/$}\kern-0.5pt
  \hbox{\phantom{}\kern-5pt$\circ$}\kern-3.0pt\hbox{$\Biggr\backslash$}
  \kern-1.5pt\hbox{\phantom{}\kern6pt{$\circ$}} }}$}}

\def\lwng{\hbox{$\vbox{\offinterlineskip{ \hbox{$\circ$}
  \kern-3.0pt\hbox{\phantom{}\kern6.0pt{$\Biggr\backslash$}}
  \kern-0.5pt\hbox{\phantom{}\kern11pt{$\circ$}}\kern-3.5pt
  \hbox{\phantom{}\kern5.0pt {$\Biggr/$}}\kern-1.0pt\hbox{$\circ$} }}$}}

\def\drwng#1#2#3{\hbox{$\vcenter{ \offinterlineskip{
  \hbox{\phantom{}\kern7pt{$\circ^{\elbl{#3}}$}}
  \kern-2.5pt\hbox{$\Biggr/$}\kern-0.5pt
  \hbox{\phantom{}\kern-5pt$\circ^{ \elbl{#1}}$}
  \kern-3.0pt\hbox{$\Biggr\backslash$}
  \kern-1.5pt\hbox{\phantom{}\kern7pt{$\circ^{\elbl{#2}}$}}  } }$}}
\def\drwngt#1#2#3{\hbox{$\vcenter{ \offinterlineskip{
  \hbox{\phantom{}\kern7pt{$\bullet^{\elbl{#3}}$}}
  \kern-2.5pt\hbox{$\Biggr/$}\kern-0.5pt
  \hbox{\phantom{}\kern-5pt$\circ^{ \elbl{#1}}$}
  \kern-3.0pt\hbox{$\Biggr\backslash$}
  \kern-1.5pt\hbox{\phantom{}\kern7pt{$\bullet^{\elbl{#2}}$}}  } }$}}

\def\dlwng#1#2#3{\hbox{$\vcenter{\offinterlineskip{ \hbox{$\lelbl{#1}\circ$}
  \kern-3.0pt\hbox{\phantom{}\kern6.0pt{$\Biggr\backslash$}}
  \kern-0.5pt\hbox{\phantom{}\kern11pt{$\lelbl{#2}\circ$}}\kern-3.5pt
  \hbox{\phantom{}\kern5.0pt
{$\Biggr/$}}\kern-1.0pt\hbox{$\lelbl{#3}\circ$}}}$}
 }

\def\rde#1#2#3{\hbox{\phantom{}\kern-4pt\hbox{$\vcenter{\offinterlineskip
\hbox
{
               \raise 4.5pt\hbox{\vrule height0.4pt width13pt depth0pt}
                \kern-1pt\vbox{ \hbox{\drwng{#1}{#2}{#3}}} }}$  }}  }
\def\rdet#1#2#3{\hbox{\phantom{}\kern-4pt\hbox{$\vcenter{\offinterlineskip
\hbox
{
               \raise 4.5pt\hbox{\vrule height0.4pt width13pt depth0pt}
                \kern-1pt\vbox{ \hbox{\drwngt{#1}{#2}{#3}}} }}$  }}  }

\def\lde#1#2#3{\hbox{$\vcenter{\offinterlineskip  \hbox{
               \dlwng{#1}{#2}{#3}\kern-4.2pt\lower0.4pt\hbox{$\vcenter{\hrule
                               width13pt}$}
               \kern-8pt\phantom{}   }}  $}}

\def\rwngb{\hbox{$\vbox{\offinterlineskip{
  \hbox{\phantom{}\kern6pt{$\bullet$}}\kern-2.5pt\hbox{$\Biggr/$}\kern-0.5pt
  \hbox{\phantom{}\kern-5pt$\bullet$}\kern-3.0pt\hbox{$\Biggr\backslash$}
  \kern-1.5pt\hbox{\phantom{}\kern6pt{$\bullet$}} }}$}}

\def\lwngb{\hbox{$\vbox{\offinterlineskip{ \hbox{$\bullet$}
  \kern-3.0pt\hbox{\phantom{}\kern6.0pt{$\Biggr\backslash$}}
  \kern-0.5pt\hbox{\phantom{}\kern11pt{$\bullet$}}\kern-3.5pt
  \hbox{\phantom{}\kern5.0pt {$\Biggr/$}}\kern-1.0pt\hbox{$\bullet$} }}$}}

\def\dbrwng#1#2#3{\hbox{$\vcenter{ \offinterlineskip{
  \hbox{\phantom{}\kern6pt{$\bullet^{\elbl{#3}}$}}
  \kern-2.5pt\hbox{$\Biggr/$}\kern-0.5pt
  \hbox{\phantom{}\kern-5pt$\bullet^{ \elbl{#1}}$}
  \kern-3.0pt\hbox{$\Biggr\backslash$}
  \kern-1.5pt\hbox{\phantom{}\kern6pt{$\bullet^{\elbl{#2}}$}}  } }$}}

\def\dblwng#1#2#3{\hbox{$\vcenter{\offinterlineskip{ \hbox{$\lelbl{#1}\bullet$}
  \kern-3.0pt\hbox{\phantom{}\kern6.0pt{$\Biggr\backslash$}}
  \kern-0.5pt\hbox{\phantom{}\kern11pt{$\lelbl{#2}\bullet$}}\kern-3.5pt
  \hbox{\phantom{}\kern5.0pt
{$\Biggr/$}}\kern-1.0pt\hbox{$\lelbl{#3}\bullet$}}}
$} }

\def\rbde#1#2#3{\hbox{\phantom{}\kern-4pt\hbox{$\vcenter{\offinterlineskip
\hbo
x{
               \raise 4.5pt\hbox{\vrule height0.4pt width13pt depth0pt}
                \kern-1pt\vbox{ \hbox{\dbrwng{#1}{#2}{#3}}} }}$  }}  }

\def\lbde#1#2#3{\hbox{$\vcenter{\offinterlineskip  \hbox{
               \dblwng{#1}{#2}{#3}\kern-4.2pt\lower0.4pt\hbox{$\vcenter{\hrule
w
idth13pt}$}
               \kern-8pt\phantom{}   }}  $}}

\def\ddgu#1.#2.{\dynk  \whtu0300{#1}\blku3000{#2}}
\def\ddgd#1.#2.{\dynk  \whtd0300{#1}\blkd3000{#2}}

\def\eddgiu#1.#2.#3.{\dynk \whtu0100{#1}\whtu1300{#2}\blku3000{#3}}
\def\eddgid#1.#2.#3.{\dynk \whtd0100{#1}\whtd1300{#2}\blkd3000{#3}}
\def\eddgiiu#1.#2.#3.{\dynk  \whtu0300{#1}\blku3100{#2}\blku1000{#3}}
\def\eddgiid#1.#2.#3.{\dynk  \whtd0300{#1}\blkd3100{#2}\blkd1000{#3}}

\def\ddfu#1.#2.#3.#4.{\dynk
\whtu0100{#1}\whtu1200{#2}\blku2100{#3}\blku1000{#4}
}
\def\ddfd#1.#2.#3.#4.{\dynk
\whtd0100{#1}\whtd1200{#2}\blkd2100{#3}\blkd1000{#4}
}

\def\eddfiu#1.#2.#3.#4.#5.{\dynk
\whtu0100{#1}\whtu1100{#2}\whtu1200{#3}\blku210
0{#4}\blku1000{#5}}
\def\eddfid#1.#2.#3.#4.#5.{\dynk
\whtd0100{#1}\whtd1100{#2}\whtd1200{#3}\blkd210
0{#4}\blkd1000{#5}}
\def\eddfiiu#1.#2.#3.#4.#5.{\dynk
\whtu0100{#1}\whtu1200{#2}\blku2100{#3}\blku11
00{#4}\blku1000{#5}}
\def\eddfiid#1.#2.#3.#4.#5.{\dynk
\whtd0100{#1}\whtd1200{#2}\blkd2100{#3}\blkd11
00{#4}\blkd1000{#5}}

\def\ddanu#1.#2.#3.#4.#5.{\dynk \whtu0100{#1}\whtu1100{#2}\whtu1100{#3}\cdots
                           \whtu1100{#4}\whtu1000{#5}}
\def\ddand#1.#2.#3.#4.#5.{\dynk \whtd0100{#1}\whtd1100{#2}\whtd1100{#3}\cdots
                           \whtd1100{#4}\whtd1000{#5}}
\def\ddandte#1.#2.#3.#4.#5.{\dynk \blkd0100{#1}\whtd1100{#2}\blkd1100{#3}\cdots
                           \blkd1100{#4}\whtd1000{#5}}
\def\ddandto#1.#2.#3.#4.#5.{\dynk \blkd0100{#1}\whtd1100{#2}\blkd1100{#3}\cdots
                           \whtd1100{#4}\blkd1000{#5}}

\def\eddanu#1.#2.#3.#4.#5.{\dynk \whtu0100{#1}\whtu1100{#2}%
                           \up1{\whtr0000{#3}}\cdots\whtu1100{#4}\whtu1000{#5}}
\def\eddand#1.#2.#3.#4.#5.{\dynk \whtd0100{#1}\whtd1100{#2}%
                           \up1{\whtr0000{#3}}\cdots\whtd1100{#4}\whtd1000{#5}}
\def\eddaid#1.#2.{\dynk\whtd0400{#1}\hskip30pt\whtd4000{#2}}

\def\eddanid#1.#2.#3.#4.#5.{\dynk \whtd0200{#1}\whtd2100{#2}%
                           \whtd1100{#3}\cdots\whtd1200{#4}\blkd2000{#5}}
\def\eddaniu#1.#2.#3.#4.#5.{\dynk \whtu0200{#1}\whtu2100{#2}%
                           \whtu1100{#3}\cdots\whtu1200{#4}\blku2000{#5}}

\def\eddaniid#1.#2.#3.#4.#5.#6.{\hbox{$\vcenter{\hbox
         {\dynk\hbox{$ \lbde{#1}{#2}{#3}\blkd1100{#4}\cdots%
          \blkd1200{#5}\whtd2000{#6} $}} }$}}
\def\eddaniiu#1.#2.#3.#4.#5.#6.{\hbox{$\vcenter{\hbox
         {\dynk\hbox{$ \lbde{#1}{#2}{#3}\blku1100{#4}\cdots%
          \blku1200{#5}\whtu2000{#6} $}} }$}}

\def\eddaiiid#1.#2.{\dynk\blkd0400{#1}\hskip30pt\whtd4000{#2}}

\def\ddbnu#1.#2.#3.#4.#5.{\dynk \whtu0100{#1}\whtu1100{#2}\whtu1100{#3}\cdots
                           \whtu1200{#4}\blku2000{#5}}
\def\ddbnd#1.#2.#3.#4.#5.{\dynk \whtd0100{#1}\whtd1100{#2}\whtd1100{#3}\cdots
                           \whtd1200{#4}\blkd2000{#5}}

\def\eddbnu#1.#2.#3.#4.#5.#6.{\dynk \lde{#1}{#2}{#3}\whtu1100{#4}\cdots
                           \whtu1200{#5}\blku2000{#6}}
\def\eddbnd#1.#2.#3.#4.#5.#6.{\dynk \lde{#1}{#2}{#3}\whtd1100{#4}\cdots
                           \whtd1200{#5}\blkd2000{#6}}

\def\ddcnu#1.#2.#3.#4.#5.{\dynk \blku0100{#1}\blku1100{#2}\blku1100{#3}\cdots
                           \blku1200{#4}\whtu2000{#5}}
\def\ddcnd#1.#2.#3.#4.#5.{\dynk \blkd0100{#1}\blkd1100{#2}\blkd1100{#3}\cdots
                           \blkd1200{#4}\whtd2000{#5}}

\def\eddcnu#1.#2.#3.#4.#5.#6.{\dynk \whtu0200{#1}\blku2100{#2}\blku1100{#3}
       \blku1100{#4}\cdots
                           \blku1200{#5}\whtu2000{#6}}
\def\eddcnd#1.#2.#3.#4.#5.{\dynk \whtd0200{#1}\blkd2100{#2}\blkd1100{#3}
       \cdots \blkd1200{#4}\whtd2000{#5}}

\def\dddnu#1.#2.#3.#4.#5.#6.{\hbox{$\vcenter{\hbox
         {\dynk\hbox{$ \whtu0100{#1}\whtu1100{#2}\cdots%
          \whtu1100{#3}\rde{#4}{#5}{#6} $}}  }$}}
\def\dddnd#1.#2.#3.#4.#5.#6.{\hbox{$\vcenter{\hbox
         {\dynk\hbox{$ \whtd0100{#1}\whtd1100{#2}\cdots%
          \whtd1100{#3}\rde{#4}{#5}{#6} $}} }$}}
\def\dddndte#1.#2.#3.#4.#5.#6.{\hbox{$\vcenter{\hbox
         {\dynk\hbox{$ \blkd0100{#1}\whtd1100{#2}\cdots%
          \blkd1100{#3}\rdet{#4}{#5}{#6} $}} }$}}
\def\dddndto#1.#2.#3.#4.#5.#6.{\hbox{$\vcenter{\hbox
         {\dynk\hbox{$ \whtd0100{#1}\blkd1100{#2}\cdots%
          \blkd1100{#3}\rdet{#4}{#5}{#6} $}} }$}}
\def\dddiv#1.#2.#3.#4.{\hbox{$\vcenter{\hbox
         {\dynk\hbox{$ \whtu0100{#1}\rde{#2}{#3}{#4}
              $}}  }$}}

\def\edddnu#1.#2.#3.#4.#5.#6.#7.#8.{\hbox{$\vcenter{\hbox
         {\dynk\hbox{$ \lde{#1}{#2}{#3}\whtu1100{#4}\cdots%
          \whtu1100{#5}\rde{#6}{#7}{#8} $}}  }$}}
\def\edddnd#1.#2.#3.#4.#5.#6.#7.#8.{\hbox{$\vcenter{\hbox
         {\dynk\hbox{$ \lde{#1}{#2}{#3}\whtd1100{#4}\cdots%
          \whtd1100{#5}\rde{#6}{#7}{#8} $}} }$}}

\def\edddniid#1.#2.#3.#4.#5.{\hbox{$\vcenter{\hbox
         {\dynk\hbox{$ \blkd0200{#1}\whtd2100{#2}\whtd1100{#3}\cdots%
          \whtd1200{#4}\blkd2000{#5} $}} }$}}
\def\edddniiu#1.#2.#3.#4.#5.{\hbox{$\vcenter{\hbox
         {\dynk\hbox{$ \blku0200{#1}\whtu2100{#2}\whtu1100{#3}\cdots%
          \whtu1200{#4}\blku2000{#5} $}} }$}}

\def\ddei#1.#2.#3.#4.#5.#6.{\hbox{$\vcenter{\hbox
       {\dynk \whtd0100{#1}\whtd1100{#3}%
       \up1{\whtr0001{#2}}\whtd1110{#4}\whtd1100{#5}\whtd1000{#6}} }$}}
\def\ddeit#1.#2.#3.#4.#5.#6.{\hbox{$\vcenter{\hbox
       {\dynk \whtd0100{#1}\blkd1100{#3}%
       \up1{\blkr0001{#2}}\whtd1110{#4}\blkd1100{#5}\whtd1000{#6}} }$}}

\def\eddei#1.#2.#3.#4.#5.#6.#7.{\hbox{$\vcenter{\hbox
       {\dynk \whtd0100{#1}\whtd1100{#3}%
       \up1{\whtr0011{#2}}\up2{\whtr0001{#7}}\whtd1110{#4}\whtd1100{#5}%
       \whtd1000{#6}} }$}}

\def\ddeii#1.#2.#3.#4.#5.#6.#7.{\hbox{$\vcenter{\hbox
       {\dynk \whtd0100{#1}\whtd1100{#3}%
       \up1{\whtr0001{#2}}\whtd1110{#4}\whtd1100{#5}\whtd1100{#6}%
       \whtd1000{#7}} }$}}
\def\ddeiit#1.#2.#3.#4.#5.#6.#7.{\hbox{$\vcenter{\hbox
       {\dynk \whtd0100{#1}\blkd1100{#3}%
       \up1{\blkr0001{#2}}\whtd1110{#4}\blkd1100{#5}\whtd1100{#6}%
       \blkd1000{#7}} }$}}

\def\eddeii#1.#2.#3.#4.#5.#6.#7.#8.{\hbox{$\vcenter{\hbox
       {\dynk \whtd0100{#8}\whtd1100{#1}\whtd1100{#3}%
       \up1{\whtr0001{#2}}\whtd1110{#4}\whtd1100{#5}\whtd1100{#6}%
       \whtd1000{#7}} }$}}

\def\ddeiii#1.#2.#3.#4.#5.#6.#7.#8.{\hbox{$\vcenter{\hbox
       {\dynk \whtd0100{#1}\whtd1100{#3}%
       \up1{\whtr0001{#2}}\whtd1110{#4}\whtd1100{#5}\whtd1100{#6}%
       \whtd1100{#7}\whtd1000{#8}} }$}}
\def\ddeiiit#1.#2.#3.#4.#5.#6.#7.#8.{\hbox{$\vcenter{\hbox
       {\dynk \whtd0100{#1}\blkd1100{#3}%
       \up1{\blkr0001{#2}}\whtd1110{#4}\blkd1100{#5}\whtd1100{#6}%
       \blkd1100{#7}\whtd1000{#8}} }$}}

\def\eddeiii#1.#2.#3.#4.#5.#6.#7.#8.#9.{\hbox{$\vcenter{\hbox
       {\dynk \whtd0100{#1}\whtd1100{#3}%
       \up1{\whtr0001{#2}}\whtd1110{#4}\whtd1100{#5}\whtd1100{#6}%
       \whtd1100{#7}\whtd1100{#8}\whtd1000{#9}} }$}}

\font\capit=cmcsc10
\font\addressit=cmcsc8
\font\eightrm=cmr8

\documentstyle[leqno,bezier,12pt]{article}

\makeatletter                                            
\renewcommand{\@begintheorem}[2]{                        
\sl \trivlist \item [\hskip \labelsep {\bf #2\ \ #1.}]   
				}                        
\makeatother                                             
\makeatletter
\def\section{\@startsection {section}{1}{\z@}{-3.5ex plus -1ex minus
 -.2ex}{1.5ex plus .2ex}{\large\bf}}

\def\subsection{\@startsection{subsection}{2}{\z@}{-3.25ex plus -1ex
minus
 -.2ex}{1.5ex plus .2ex}{\normalsize\it}}

\newcommand{\numberequationsassubsubsections}

\newtheorem{prop}{Proposition}[section]
\newtheorem{lemm}[prop]{Lemma}
\newtheorem{mumlemm}[prop]{Atiyah-Mumford lemma}
\newtheorem{theo}[prop]{Theorem}
\newtheorem{cor}[prop]{Corollary}
\newtheorem{conj}[prop]{Conjecture}
\newtheorem{rem}[prop]{\it Remark}
\newtheorem{rems}[prop]{\it Remarks}
\newtheorem{ex}[prop]{Example}

\begin{document}

\title{Spin Verlinde spaces and Prym theta functions}
\author{W.M. Oxbury}
\date{}

\maketitle

\centerline{\it Department of Mathematical Sciences}
\centerline{\it University of Durham}
\bigskip

\section{Introduction}

The Verlinde formula is a remarkable---and potentially very useful---new tool
in the geometry of algebraic curves which is borrowed from conformal field
theory. In the first instance it is a trigonometric expression which assigns a
natural number $N_l(G,g)$ to data consisting of a semisimple algebraic group
$G$, a nonnegative integer $g$ and an auxiliary integer $l \in \z$. In physics
$N_l(G,g)$ is interpreted as the dimension of the space of conformal blocks at
level $l$ in the Wess-Zumino-Witten model of conformal field theory on a
compact Riemann surface of genus $g$. In algebraic geometry this `space of
conformal blocks' is identified with a vector space
of the form $H^0(\mm_C(G),\ll^l)$, where $\mm_C(G)$ is the moduli scheme (or
stack) of semistable principal $G$-bundles over $C$, and $\ll$ is an ample line
bundle on this moduli space generating the Picard group.

The feature of the Verlinde formula which motivates this paper
is its `numerology'. Namely, when one computes the numbers $N_l(G,g)$ for the
classical simple groups one finds that they obey interesting identities which
lead one to certain conjectures about the geometry of $\mm_C(G)$. (See
\cite{OW}.) In this article we are concerned with identities linking the
complex spin groups $G=\spin_m$ with the configuration of principally polarised
Prym varieties associated to the curve $C$ via its unramified double covers.
This connection was first observed in \cite{O} for the odd spin groups
$\spin_{2n+1}$; here we shall give a systematic account of both the odd and
even cases.

The moduli space $\mm_C(SO_m)$ has, for $m\geq 3$, two connected components
labelled by the second Stiefel-Whitney class $w_2$ of the bundles. Each of
these components has a $J_2(C)$-Galois cover which is a moduli space for
Clifford bundles with fixed spinor norm line bundle $\xi \in \pic(C)$, and
whose isomorphism class depends only on $\deg \xi $ mod 2. When $\xi = \oo_C$
the Galois cover is precisely $\mm_C(\spin_m)$ since by definition $\spin_m$ is
the kernel of the spinor norm; the `sister' space which arises when $\deg \xi$
is odd, we denote by $\mm_C^-(\spin_m)$. In each case the fibre $J_2(C)$ over
$\mm_C(SO_m)$ parametrises liftings of a given $SO_m$-bundle to a Clifford
bundle with given spinor norm $\xi$ with $\deg \xi \equiv w_2$ mod 2.

Note that for low values of $m$ we recover well-known moduli spaces of vector
bundles: for example $\mm_C(\spin_3) = \su(2,0)$ and $\mm_C^-(\spin_3) =
\su(2,1)$ (where $\su(n,d)$ is the moduli space of rank $n$ vector bundles with
fixed determinant of degree $d$); while $\mm_C(\spin_6) = \su(4,0)$ and
$\mm_C^-(\spin_6) = \su(4,2)$.

We shall consider the theta line bundle $\Theta(\c^m)$ on these varieties,
coming from the standard orthogonal representation $\c^m$ of the Clifford
group, and use the Verlinde formula to count its sections.

These preliminary ideas---the spin moduli spaces, theta line bundles on them
and the Verlinde calculations---occupy the first five sections of the paper.
The central observations on which the paper is based is contained in section
\ref{numerology}.
This is that the dimension of $H^0(\mm_C(\spin_m),\Theta(\c^m))$ coincides with
that of the direct sum of the spaces of {\it even} level $m$ theta functions on
all the Prym varieties (including the Jacobian itself); while the dimension of
$H^0(\mm_C^-(\spin_m),\Theta(\c^m))$ is equal to that of the direct sum of the
spaces of {\it odd} level $m$ theta functions on the Prym varieties. The
precise statement (see theorem \ref{numer} and table (\ref{dims}))
depends on the parity of $m$: when $m$ is even we have to take theta functions
not just on the Prym varieties $\peta$, $\eta \in J_2(C)\backslash \{\oo\}$,
but on the two component abelian subvarieties
$$
\peta \cup \peta^- = \nm^{-1}K_C \subset J^{2g-2}(\ctil),
$$
where $\ctil \rightarrow C$ is the double cover corresponding to each $\eta$.

There are two remarks to make about this feature of the even spin groups. The
first is that it
is natural in the sense that, whereas level $m$ theta functions are
well-defined on $\peta$ for all $m$---there is a canonically defined theta
divisor $\xieta$---on $\peta^-$ they are well-defined only if $m$ is even.
The second is that
it corresponds to a direct sum decomposition of
$H^0(\mm_C(\spin_m),\Theta(\c^m))$ into two pieces when $m$ is
even---apparently the eigenspaces under the involution of the moduli space
$\mm_C(\spin_m)$ corresponding to reflection of the Dynkin diagram. (See
section \ref{dynkin}.)

In the remaining sections of the paper we make sense of the observations of
section \ref{numerology} by constructing homomorphisms:
$$
H^0(\mm^{\pm}(\spin_{2n+1}), \Theta(\c^{2n+1}))^{\vee} \rightarrow
\sum_{\eta \in J_2}H^0_{\pm}(\peta,({2n+1})\xieta)
$$
$$
\begin{array}{rcl}
H^0(\mm^{\pm}(\spin_{2n}), \Theta(\c^{2n}))^{\vee} &\rightarrow&
\displaystyle
\sum_{\eta \in J_2}H^0_{\pm}(\peta,2n\xieta) \\
&&\ \ \ \ \ \ \displaystyle
\oplus \sum_{\eta \not= 0}H^0_{\pm}(\peta^-,2n\xieta).\\
\end{array}
$$
Table (\ref{dims}) asserts that each of these maps is {\it between vector
spaces of equal dimension}. It is necessary to emphasise that the validity of
the right-hand column of the table is dependent on a natural-seeming conjecture
\ref{conj} for a Verlinde formula on $\mm^-(\spin_m)$. When $m=3$ this is the
Verlinde formula for rank 2 vector bundles of odd degree due to Thaddeus
\cite{T}; further evidence for this `twisted' formula is given in \cite{OW}.

It is therefore natural to expect that the above maps
are isomorphisms; for $m\geq 5$ this is not known, though we hope to return to
the question in a later paper.
Cases of low $m$, on the other hand, where the spin moduli spaces can be
identified with more familiar moduli spaces of vector bundles, are examined
individually in section \ref{numerology}.

Finally, when $m$ is odd it is known (see \cite{LS}) that $\Theta(\c^m) = 2\pp$
where the `Pfaffian' line bundle $\pp$ generates the Picard group. In section
\ref{pfaffian} we observe---though this remark is independent of the rest of
the paper---that the space of sections of
this line bundle has a basis labelled by the even theta characteristics of the
curve, which directly generalises that constructed by Beauville in \cite{B2}
for the case $m=3$.

\medskip
\noindent
{\it Acknowledgements:} In writing this paper the author has benefited
greatly from conversations with B. van Geemen, C. Pauly, S. Ramanan and C.
Sorger, to all of whom he expresses his gratitude.

\section{Moduli spaces of principal bundles on a curve}
\label{mod}

In this section we shall give a brief account of the moduli spaces of
semistable principal bundles over a curve, following [R1], [DN], [KNR].

We begin with a smooth projective complex curve $C$ of genus $g \geq 2$, and a
complex connected reductive algebraic group $G$; and we consider algebraic
principal $G$-bundles $E\rightarrow C$. Topologically such bundles are
classified by the fundamental group of $G$.

Just as for vector bundles, one has notions of stability, semistability and
S-equivalence for algebraic $G$-bundles, and for stable bundles S-equivalence
is the same as isomorphism. (We shall recall in a moment the definition of
stability, but it will not be necessary here to define S-equivalence.)
The basic result of Ramanathan [R2] is then the following.

\begin{theo}
Given $C,G$ as above and an element $\gamma \in \pi_1(G)$, there exists a
normal irreducible projective variety $\mm(G,\gamma)$ which is a coarse moduli
space for families of semistable $G$-bundles of type $\gamma$
on $C$, modulo S-equivalence.
\end{theo}

Moreover, one has
$$
\dim \mm(G,\gamma) = (g-1)\dim G + \dim Z(G),
$$
and $\mm(G,\gamma)$ is unirational when $G$ is a simple group [KNR].

The basic construction with principal bundles is the following. If $E$ is a
$G$-bundle, and $\rho:G\rightarrow {\rm Aut}(X)$ any left $G$-space, then we
can form a bundle $E(X) = E\times_{\rho} X$ with fibre $X$. In case $X=G/P$ is
a homogeneous coset space, a section $\sigma: C \rightarrow E(G/P)$ is called a
reduction of the strucure group of the bundle to the subgroup $P$. When
$P\subset G$ is a maximal parabolic, $E(G/P)\rightarrow C$ can be thought of as
a `generalised Grassmannian bundle'. Then by definition, $E$ is {\it
semistable} if and only if
$$
\deg \sigma^* T^{\rm vert}E(G/P) \geq 0
\qquad
\hbox{\sl for all maximal parabolics $P\subset G$,}
$$
where $T^{\rm vert}$ denotes the vertical tangent bundle.

On the other hand, if $\pi : G' \rightarrow G$ is a group epimorphism then we
can view $X=G$ as a left $G'$-space via $\pi$, and so form a $G$-bundle
$E=F(G)$ from any $G'$-bundle $F$. $F$ is said to be a {\it lift} of $E$. In
particular, if $G'$ is a central extension of $G$ then there is a bijection
between maximal parabolics $P\subset G$ and maximal parabolics $P'=\pi^{-1}
P\subset G'$, and moreover $F(G'/P') \cong E(G/P)$ if $F$ is any lift of $E$.
Consequently:

\begin{lemm}
\label{1.2}
If $E$ is a $G$-bundle and $F$ a lift of $E$ to a central extension of $G$ then
$E$ is stable (resp. semistable) if and only if $F$ is.
\end{lemm}

Finally, of course, we can take for the $G$-space $X$ a finite-dimensional
representation $\rho:G\rightarrow GL(V)$, to obtain a vector bundle $E(V)$.
In the case when $G = GL_n$ and $V = \c^n$ is the standard representation, the
notions of stability, semistability and S-equivalence are
the same for the principal bundle $E$ as for the vector bundle $E(V)$. Thus
we shall write $\uu(n,d) = \mm(GL_n,d)$, for $d\in \pi_1(GL_n) \cong \z$; this
is the moduli space of semistable vector bundles of rank $n$ and degree $d$.

Consider now the determinant morphism
$$
\det : \uu(n,d) \rightarrow J^d(C).
$$
This is a fibration and we shall, as is usual, denote the isomorphism class of
the fibre by $\ss\uu(n,d) = \ss\uu_C(n,d) $; or $\ss\uu_C(n)$ when $d=0$.

One knows from [DN] that, via det, $\uu(n,d)$ has Picard group
$$
\pic \ \uu(n,d) \cong \pic \ J^d(C) \oplus \z\{\Theta_{n,d}\}
$$
where $\Theta_{n,d}$ is an ample line bundle on the fibres
constructed as follows. It will be convenient, to begin with, to assume that
$n$ divides $d$, i.e. that we are dealing with vector bundles of integral
slope.

Consider first an arbitrary family $F \rightarrow C\times S$ of semistable
vector bundles on $C$ with rank $n$, degree $d$ and slope $\mu = d/n \in \z$,
as above; and we construct a line bundle $\Theta (F) \rightarrow S$, functorial
with respect to base change $S' \rightarrow S$, in the following way. Let $\pi
: C\times S \rightarrow S$ be the projection. Then (at least Zariski locally)
there is a homomorphism of locally free sheaves on $S$, $\phi :K^0 \rightarrow
K^1$, having the direct images of $F$ under $\pi$ as kernel and cokernel:
\begin{equation}
\label{1.3}
0\rightarrow R^0_{\pi}F
\rightarrow K^0\map{\phi} K^1\rightarrow R^1_{\pi}F\rightarrow 0.
\end{equation}
Moreover, the determinant line bundle
$$
\textstyle
{\rm Det}(F) = (\bigwedge^{\rm top} K^0)\dual \otimes (\bigwedge^{\rm top} K^1)
$$
is well-defined and functorial with respect to base change. If $\mu = g-1$ we
write $\Theta(F) = {\rm Det}(F)$; and this has a canonical section $\det \phi$,
so that in this case $\Theta(F)$ is represented  by a canonical Cartier divisor
on $S$.
Otherwise $\Theta(F)$ is defined to be a suitable twist of ${\rm Det}(F)$ such
that
$$
\Theta(F\otimes \pi^*L) = \Theta(F)
$$
for any line bundle $L\rightarrow S$; i.e. $\Theta$ respects equivalence of
families.

Now in the case $\mu = g-1$ it is shown in [DN] that the functor $\Theta$ is
in fact
represented in moduli space
by a global Cartier divisor
$$
\Theta_{n,n(g-1)} = {\rm Closure}\{{\rm stable}\ V| H^0(V) = H^1(V) \not= 0 \}
\subset \uu(n,n(g-1)).
$$
In other words $\Theta(F) = f^* \Theta_{n,n(g-1)}$ where $f: S \rightarrow
\uu(n,n(g-1))$ is given by the coarse moduli property.

For the general case ($\mu \in \z$ still) one chooses a line bundle $L\in
\pic(C)$ with degree chosen so that we get a morphism
$$
\uu(n,d) \map{\otimes L} \uu(n,n(g-1)).
$$
Now set
$$
\Theta_L = (\otimes L)^*\Theta_{n,n(g-1)}.
$$
The dependence of $\Theta_L$ on $L$ is then given by (\ref{1.5}) below, which
is a consequence of:

\begin{lemm}
\label{1.4}
View $J^0(C)$ as a subgroup of $\pic(J^d(C))$ by $L \mapsto \Phi^{-1}\otimes
T^*_L\Phi$, where $\Phi$ is any line bundle representing the principal
polarisation on $J^d(C)$. Then for any family $F\rightarrow C\times S$ as
above, and any $L\in J^0(C)$, we have
$$
\Theta(F\otimes pr_C^*L) = {\det}^*(L) \otimes \Theta(F)
$$
where $\det : S\rightarrow \uu(n,d) \rightarrow J^d(C)$.
\end{lemm}

It follows easily from this that when $L,L'$ have the same degree,
$\Theta_L$ and $\Theta_{L'}$ are related by:
\begin{equation}
\label{1.5}
\Theta_{L'} = {\det} ^*(L'L^{-1}) \otimes \Theta_L \in \pic\ \uu(n,d).
\end{equation}
We now set $\Theta_{n,d} = \Theta _L$: if $d\not= n(g-1) $ this depends on $L$,
{\it but by (\ref{1.5}) its restriction to the fibres of $\det :\uu(n,d)
\rightarrow J^d(C)$
is independent of $L$}.

\medskip

In order to consider bundles of general degree, i.e. non-integral slope, it
is necessary to twist by bundles $L$ of higher rank:
$$
\uu(n,d) \map{\otimes L} \uu(rn,rn(g-1)),
$$
where $L$ has rank $r$. It is easy to check that the necessary and sufficient
condition for arranging slope $g-1$ on the right is that:
\begin{equation}
\label{minr}
r\in {n\over {\rm gcd}(n,d)} \z
\end{equation}
The line bundle
$\Theta_L$ may now be defined in the same way as above. In this more general
situation (\ref{1.5}) becomes:
\begin{equation}
\label{1.5'}
\Theta_{L'} = {\det} ^*(\det L' \otimes \det L\dual) \otimes \Theta_L \in \pic\
\uu(n,d),
\end{equation}
Consequently the restriction of
$\Theta _L$ to the fibres of $\det :\uu(n,d) \rightarrow J^d(C)$
is again independent of the choice of $L$ with given rank $r$; and we set
$\Theta_{n,d} = \Theta _L$ for any $L$ with $r = {n/ {\rm gcd}(n,d)}$. This is
the required generator of the Picard group.

Note that if in this construction $L,L'$ are two vector bundles of different
ranks $r<r'$ (both satisfying (\ref{minr}), and the degrees of $L$ and $L'$
chosen suitably) then
\begin{equation}
\label{difr}
\Theta_{L'} = \Theta_L^{\otimes r'/r}.
\end{equation}

\medskip

Finally, suppose that we are given a family $E\rightarrow C\times S$ of
semistable $G$-bundles, and a representation $\rho:G\rightarrow SL(V)$, where
$\dim V =n$.
We shall suppose that $\rho$ satisfies the condition:
\begin{equation}
\label{ss}
Z(G)_0 \subset \ker \rho
\end{equation}
We can form the family of vector bundles $E(V) \rightarrow C\times S$;
and by \cite{R2} proposition 2.17 the condition (\ref{ss}) guarantees that
these vector bundles are semistable.

We thus obtain a theta line bundle $\Theta(E(V))\rightarrow S$, and since
$E(V)$ has trivial determinant on the fibres of $\pi:C\times S\rightarrow S$ we
deduce from lemma \ref{1.4} the following corollary, which will be needed
later:

\begin{cor}
\label{1.6}
 For $E\rightarrow C\times S$ and  $\rho:G\rightarrow SL(V)$ as above, and for
any $L\in J^0(C)$ one has
$$
\Theta(E(V)\otimes pr_C^*L) = \Theta(E(V)).
$$
\end{cor}

Globally $\rho$ satisfying (\ref{ss}) induces a morphism
$$
\rho_* : \mm(G,\gamma) \rightarrow \mm(SL_n) \hookrightarrow \uu(n,0)
$$
and the functor $E\mapsto \Theta(E(V))$ is represented by the (well-defined)
line bundle
$$
\Theta(V) := (\rho_*)^* \Theta_{n,0} \in \pic \ \mm(G,\gamma).
$$

Note that if we let $j: \mm(SL_n) \hookrightarrow \uu(n,0)$ denote the
inclusion which identifies $\mm(SL_n)$
with the moduli space of vector bundles of rank $n$ and trivial determinant,
via the standard representation $\c^n$, then by construction
$$
\Theta(\c^n) = j^*\Theta_{n,0}.
$$

\section{Clifford bundles}
\label{clif}

Let us consider again the fibration
$$
\det :\ \uu_C(n,d) \rightarrow J^d(C);
$$
induced, that is, by
the determinant homomorphism $GL_n \rightarrow \c^*$. The fibres of these maps
are, up to isomorphism, the $n$ moduli varieties $\ss\uu_C(n,d)$, for $d\in \z
/n$.
In this section we shall describe an alternative generalisation of this
situation for $n=2$, obtained by replacing $GL_2$ not by $GL_n$, but
by the special Clifford group of a nondegenerate quadratic form. First we need
to recall some basic Clifford theory.

\subsection{The special Clifford group}

Let $Q$ be a nondegenerate quadratic form on a complex vector space $V$ of
finite dimension $m$;
let $A=A(Q)$ be its Clifford algebra and $A^+$ the even Clifford algebra.
Recall that these
can be expressed as matrix algebras as follows.

If $m=2n$ is even then for any $n$-dimensional isotropic subspace $U\subset V$
one has
\begin{equation}
\label{evenA}
\textstyle
A \cong \endo \bigwedge U;
\qquad
A^+ \cong \endo \bigl( \bigwedge^{\rm even} U \bigr) \oplus \endo \bigl(
\bigwedge^{\rm odd} U \bigr).
\end{equation}

If, on the other hand, $m=2n+1$ is odd then for any direct sum decomposition
$V=U \oplus U'\oplus \c$ where $U,U'$ are $n$-dimensional isotropic subspaces
one has

\begin{equation}
\label{oddA}
\textstyle
A \cong \endo \bigwedge U \oplus \endo \bigwedge U';
\qquad
A^+ \cong \endo \bigwedge U .
\end{equation}
The `principal involution' of $A$ is $\alpha : x\mapsto -x$ for $x\in V$, i.e.
is $\pm 1$ on $A^{\pm}$ respectively. The `principal anti-involution' $\beta$
is the identity on $V$ and reverses the direction of multiplication:
$\beta(x_1\ldots x_r) = x_r\ldots x_1$. Then the {\it Clifford group} is
$$
 C (Q) = \{s\in  A^{*} | \alpha (s)Vs^{-1} \subset V\},
$$
where $A^{*}\subset A$ denotes the group of units; and the {\it special
Clifford group} is
$$
SC(Q) = C(Q) \cap A^+.
$$
For $s\in C(Q)$ the transformation $\pi_s: x\mapsto \alpha(s)xs^{-1}$ of $V$ is
orthogonal---this is because $C(Q)$ is generated by $x\in V\cap C(Q)$, for
which $\pi_x$ is just minus the reflection in the hyperplane $x^{\perp}$. Thus
one has a group homomorphism $\pi:C(Q) \rightarrow O(Q)$, which
has the following properties.

\begin{prop}
\label{2.1}
\begin{enumerate}
\item $\ker \pi = \c^*$;
\item $\pi(C(Q)) = O(Q)$ and $\pi(SC(Q)) = SO(Q)$.
\end{enumerate}
\end{prop}


\begin{cor}
\label{2.2}
$SC(Q)$ is a connected reductive algebraic group.
\end{cor}

The {\it spinor norm} is the group homomorphism
$$
\begin{array}{rcl}
\nm \ :\ SC(Q) &\rightarrow& \c^{*} \\
         s &\mapsto& \beta(s)s.\\
\end{array}
$$
Equivalently $\nm (x_1\ldots x_r) =  Q(x_1)\cdots Q(x_r)$ for $x_1,\ldots ,x_r
\in V$. Then by definition $\spin (Q)= \ker \nm$.

Note that multiplication by scalars induces a double cover
\begin{equation}
\label{rescaling}
\{\pm(1,1)\} \rightarrow \c^* \times \spin(Q) \rightarrow SC(Q).
\end{equation}

{}From now on we shall write $C_m$, $SC_m$, $\spin_m$ instead of
$C(Q)$, $SC(Q)$, $\spin(Q)$ when $Q$ is the standard quadratic form on $\c^m$.
Then
(using (\ref{2.1})) one has the following commutative diagram of short exact
sequences:

\begin{equation}
\label{2.3}
\matrix{
1&\rightarrow&\c^{*}&\rightarrow&SC_m&\map{\pi}&SO_m&\rightarrow&1\cr
&&&&&&&&\cr
&&\uparrow&&\uparrow&&\Vert&&\cr
&&&&&&&&\cr
1&\rightarrow&\z/2&\rightarrow&\spin_m&\rightarrow&SO_m&\rightarrow&1.\cr
}
\end{equation}

\begin{prop}
\label{2.4}
For $m\geq 3$:
\begin{enumerate}
\item $SC_m$ has centre
$
Z(SC_m) =
\cases{\c^{*} & if $m$ is odd \cr
       \c^{*} \times \z/2 & if $m$ is even;\cr}
$
\item $SC_m$ has fundamental group $\pi_1(SC_m) = \z$; and this maps
isomorphically to $\pi_1(\c^*) = \z$ under the spinor norm.
\end{enumerate}
\end{prop}

\pf (i) If $m$ is odd the centre of $SC_m$ must be contained in---and hence
equal to---the kernel $\c^*$ of the surjection onto $SO_m$, since the latter
has trivial centre.

If, on the other hand,
 $m$ is even, then by the same token $Z(SC_m)$ is contained in
$\pi^{-1}(Z(SO_m))$ where $\pi$ denotes the surjection to $SO_m$. In this case
$SO_m$ has centre $\{\pm 1\}$. As before everything in $\pi^{-1}(1) = \c^*$
is central; while if $\{ e_1, \ldots , e_m \} \subset \c^m$ is any orthonormal
basis then the product $e_1 \ldots  e_m \in SC_m$ spans $\pi^{-1}(-1) \cong
\c^*$. Since $m$ is even this product anticommutes with each $e_i$, and
therefore {\it commutes} with all elements of $A^+$. So $\pi^{-1}\{\pm 1\}
\cong \c^{*} \times \z/2$ is contained in and therefore equal to the centre.

(ii) From the exact homotopy sequence of the fibration in the upper sequence of
(\ref{2.3}), and the vanishing of $\pi_2(\c^*)$, we have a non-split extension
$$
\ses{\z}{\pi_1(SC_m)}{\z /2}.
$$
Since the fundamental group of a Lie group is abelian it follows that the only
possibility is $\pi_1(SC_m) = \z$. The last part now follows from the fact that
$\spin_m$ is simply-connected.
\qed

\subsection{The spin moduli spaces}

{}From proposition \ref{2.4} part 2 we see that there is for each $d\in \z =
\pi_1(SC_m)$ a morphism
induced by the spinor norm, and which we shall denote in the same way:
$$
\nm : \mm (SC_m, d) \rightarrow J^d(C).
$$
Moreover, when $m=3$ this is nothing but the determinant morphism for rank 2
vector bundles (see example \ref{clifm=3} below).
And just as for rank 2 vector bundles, one has:

\begin{prop}
\label{2.5}
For $d\in \z$ and $L\in J^d(C)$, the isomorphism class of the subscheme $\nm
^{-1}(L) \subset \mm(SC_m,d)$ depends only on $d$ mod 2.
\end{prop}

Actually this is essentially trivial and is proved in the same way as for rank
2 vector bundles, once one observes that multiplication of $SC_m$ by its centre
(proposition \ref{2.4}) induces a natural generalisation of the
tensor product operation of Clifford bundles by line bundles if $m$ is odd, and
if $m$ is even by {\it pairs} $(N,\eta)$ where $N$ is a line bundle and $\eta
\in H^1(C,\z/2) = J_2(C)$. (See also \cite{R}.)
We shall write, respectively, $N\otimes E$ and $(N,\eta)\otimes E$ for this
product.
It follows from the definition of the spinor norm that
$$
\nm(N\otimes E) = N^2 \otimes \nm (E),
\quad {\rm resp.}
\quad \nm((N,\eta)\otimes E) = N^2 \otimes \nm (E).
$$

So to prove \ref{2.5}: suppose $d=\deg L \equiv d'=\deg L'$ mod 2, and write
$L' = N^{2k}\otimes L$ for some $N\in \pic^1(C)$, $k\in \z$. Then the map
$\mm(SC_m,d)
\rightarrow \mm(SC_m,d')$ given by $E\mapsto N\otimes E$ (resp. $(N,\oo)\otimes
E$) restricts to an isomorphism $\nm^{-1}(L) \cong \nm^{-1}(L')$.
\qed

We shall therefore introduce the notation:
$$
\begin{array}{rcl}
\mm^+(\spin_m) &=& \mm(\spin_m) = \nm^{-1}(\oo_C);\\
\mm^-(\spin_m) &=& \nm^{-1}(\oo_C(p)),\\
\end{array}
$$
where $p\in C$ is any point of the curve.

\begin{rem}\rm
\label{w2}

The group $SO_m$ has fundamental group $\z/2$; so the moduli space of
semistable $SO_m$-bundles has two irreducible components $\mm(SO_m,0) =
\mm^+(SO_m)$ and $\mm(SO_m,1) = \mm^-(SO_m)$ distinguished by the second
Stiefel-Whitney class $w_2$.
It is not hard to show (see \cite{O}) that $\deg \nm(E) \equiv w_2(E)$ mod 2
and that
$\mm^{\pm}(\spin_m)$ are naturally Galois covers of these components:
$$
\mm^{\pm}(\spin_m) \map{J_2(C)} \mm^{\pm}(SO_m).
$$
Note that these maps respect stability---this is a special case of lemma
\ref{1.2}.
\end{rem}

In later sections we shall be interested in the theta line bundle
$\Theta(\c^m)$ associated to the orthogonal representation of $SC_m$. Note that
by \ref{2.1} and \ref{2.4} the condition (\ref{ss}) is satisfied for $m\geq 3$,
so that $\Theta(\c^m)$ is defined everywhere on $\mm(SC_m)$ and hence on both
$\mm^{\pm}(\spin_m)$.

\begin{ex}
\label{clifm=2} $\bf m = 2.$ \rm
Here $SC_2 \cong \c^* \times \c^*$, the spinor norm is $\nm : (a,b) \mapsto ab$
and the orthogonal representation $SC_2 \rightarrow SO_2 \cong \c^*$ is $(a,b)
\mapsto a/b$.
It follows that $\mm(SC_2) \cong \pic(C)\times \pic(C)$ and each of
$\mm^{\pm}(\spin_2) \cong \pic(C)$.

$\Theta(\c^2)$ is by definition obtained by pulling back $\Theta_{2,0}$ under
$$
\rho_*: \mm(SC_2) = \pic(C)\times \pic(C) \rightarrow \mm(SL_2).
$$
But this map sends a pair of line bundles $(L,N)$ to the vector bundle
$LN^{-1}\oplus L^{-1}N$, and this is semistable only if $\deg L = \deg N$.
(Note that condition (\ref{ss}) fails for this case!) It follows that the
morphism $\rho_*$, and hence
$\Theta(\c^2)$, is defined only on $\pic(C)\times_{\deg} \pic(C)$; and
$\Theta(\c^2)$ is thus defined only on the degree 0 component $J(C)$ of
$\mm(\spin_2) = \pic(C)$ and is not defined on $\mm^{-}(\spin_2)$.

On the other hand, $\mm(\spin_2) \rightarrow \mm(SO_2)$ is the squaring map
$[2]$ on line bundles, and so one sees that on $J(C)$ the orthogonal theta
bundle is:
$$
\Theta(\c^2) = [2]^*(2\theta) = 8 \theta,
$$
where $\theta$ is the theta divisor on $J(C)$.
\end{ex}

\begin{ex}
\label{clifm=3} $\bf m = 3.$ \rm
This is in many ways the most important case. It is well-known that $SC_3$ is,
via (\ref{oddA}), equal to the group of units $GL_2 \subset A^+$, and that the
spinor norm is the determinant homomorphism so that $\spin_3 = SL_2$. (The
representation on $SO_3$ is then precisely the action of $GL_2$ on $S^2 \cong
\bp_1$ by M\"obius transformations, via stereographic projection.) Thus
$\mm(\spin_3) \cong \su(2)$
and $\mm^-(\spin_3) \cong \su(2,1)$.

Let $\ll_d = \Theta_{2,d}$ be the (ample) generators of the Picard groups of
these varieties, for $d=0,1$ respectively, as described in the previous
section. From the discussion there we can view
$\ll_1$ ( $=\Theta_{L'}$ for a suitable rank 2 vector bundle $L'$) as defined
on all of $\mm(GL_2)$ while $\ll_0$ ($= \Theta_L$ for a suitable line bundle
$L$) is a line bundle defined only on the components of even degree; and by
(\ref{difr}) $\ll_1 = \ll_0^2$ on these components.

The orthogonal theta bundle is then $\Theta(\c^3) = \ll_1^2 = \ll_0^4$---see
(\ref{di2}) below.
\end{ex}

\begin{ex}
\label{clifm=4} $\bf m = 4.$ \rm
Again it is very well-known that $\spin_4 = SL_2 \times SL_2$. Via
(\ref{evenA}) the special Clifford group $SC_4 \subset GL_2 \times GL_2$ is the
subgroup consisting of pairs of matrices $(A,B)$ such that $\det A = \det B$;
the spinor norm is then the common $2\times 2$ determinant. Thus $\mm(SC_4) =
\mm(GL_2) \times_{\det} \mm(GL_2)$ while $\mm^{\pm}(\spin_4)$ are $\su(2)\times
\su(2)$ and $\su(2,1)\times \su(2,1)$ respectively.

The orthogonal representation induces $\mm(SC_4) \rightarrow \mm(SL_4)$ mapping
a pair of rank 2 vector bundles $(E,F)$ to $E\otimes F^*$, and from this it
follows that the orthogonal theta bundle is (with the notation of the previous
example)
$$
\Theta (\c^4) = pr_+^* \ll_1 \otimes pr_-^* \ll_1 = pr_+^* \ll_0^2 \otimes
pr_-^* \ll_0^2,
$$
where $pr_{\pm}$ denote the respective projections.

\end{ex}

\begin{ex}
\label{clifm=6} $\bf m = 6.$ \rm
In this case $\spin_6 \cong SL_4$, and one may show that the subgroup (using
(\ref{evenA}) once again) $SC_6 \subset GL_4 \times GL_4$ is the image of the
homomorphism
$$
\begin{array}{rcl}
\c^* \times SL_4 &\rightarrow& GL_4 \times GL_4\\
(\lambda,A) &\mapsto& (\lambda A,\lambda A^{adj,t})\\
\end{array}
$$
where $A^{adj,t}$ is the matrix of signed cofactors and $A^{adj} = \det A
\times A^{-1} = A^{-1}$. (Compare with (\ref{rescaling}).)

Note that projection to the first factor---the first half-spinor
representation---induces a double cover
$
\{(1,\pm 1)\} \subset SC_6 \map{pr_1} GL_4.
$
Moreover, the lift of the determinant function to this double cover has a
square root---namely, the spinor norm
$\nm : (\lambda A,\lambda A^{adj,t}) \mapsto \lambda^2$, where
$\nm(a)^2 = \det pr_1(a)$ for $a\in SC_6$.

At the level of bundles this says we have a commutative diagram
$$
\matrix{
\mm(SC_6)&\map{pr_1}&\mm(GL_4)\cr
&&\cr
\nm\downarrow&&\downarrow\det\cr
&&\cr
\pic(C)&\map{[2]}&\pic(C)\cr}
$$
where each of the horizontal maps has fibre $J_2(C)$. In particular one sees
that $\mm(\spin_6) \cong \su(4)$
and $\mm^-(\spin_6) \cong \su(4,2)$.

As in example \ref{clifm=3}, let $\ll_d = \Theta_{2,d}$ be the ample generators
of the Picard groups of these varieties, for $d=0,2$ respectively.  Then
$\ll_2$ ( $=\Theta_{L'}$ for a suitable rank 2 vector bundle $L'$) is defined
on all components of $\mm(GL_2)$ of even degree,
while $\ll_0$ ($= \Theta_L$ for a suitable line bundle $L$) is defined only on
the components of degree $\equiv 0$ mod 4; and by (\ref{difr}) $\ll_2 =
\ll_0^2$ on these components. Using (\ref{di2}) from section \ref{vformula}
below,
the orthogonal theta bundle is then $\Theta(\c^6) = \ll_2 = \ll_0^2$.

\end{ex}

\section{Orthogonal bundles and theta characteristics}

In this section we shall gather together various properties of orthogonal
bundles, some possibly well-known, which will be needed later on.

\subsection{Isotropic line subbundles}

To begin, consider any $SO_m$-bundle $E$ and its associated orthogonal vector
bundle $E(\c^m)$.
Recall that stability of $E$ is equivalent to the condition that $\deg F <0$
for all {\it isotropic} vector subbundles $F\subset E(\c^m)$. This holds, in
particular, if $E(\c^m)$ is stable as a vector bundle.

On the other hand, for
{\it any} subbundle $F\subset E(\c^m)$, the direct sum
with its orthogonal complement fits into an exact sequence
$$
\ses{N}{F\oplus F^{\perp}}{M}
$$
where $N,M$ are the subbundles generically generated by $F\cap F^{\perp}$ and
$F+ F^{\perp}$ respectively. In particular, when $N=0$ this gives rise to an
orthogonal splitting $E(\c^m) = F\oplus F^{\perp}$.
Applying this idea inductively Ramanan shows (\cite{R} proposition 4.5):

\begin{lemm}
\label{ramanan}
$E$ is a stable $SO_m$-bundle if and only if $E(\c^m)$ is an orthogonal direct
sum $E(\c^m)= F_1 \oplus \cdots \oplus F_k$ where the $F_i$ are pairwise
non-isomorphic stable vector bundles.

In particular $E(\c^m)$ is a stable vector bundle for generic $E\in \mm(SO_m)$.
\end{lemm}

\noindent
(The last assertion here follows from a simple dimension count.)

Now suppose that $F\subset E(\c^m)$ is a {\it line} subbundle. Then $N=0$
precisely when $F$ is non-isotropic, and thus $E(\c^m)$ splits in this case,
and so fails to be stable as a vector bundle. This shows:

\begin{cor}
\label{alliso}
For generic $E\in \mm(SO_m)$ every line subbundle $F\subset E(\c^m)$ is
isotropic.
\end{cor}

In the next lemma we observe that the bound $\deg F <0$ for isotropic
subbundles of $E(\c^m)$---again restricting to line subbundles---can in fact be
improved generically:

\begin{lemm}
\label{nagata}
For generic stable $E\in \mm^{\pm}(SO_m)$ every isotropic line subbundle
$F\subset E(\c^m)$ satisfies
$
\deg F \leq -g+1.
$
\end{lemm}

\pf
We simply count dimensions of those bundles $E\in \mm(SO_m)$ which can possess
an isotropic line subbundle $F\subset E(\c^m)$ with $\deg F = d$. Let
$F^{\perp}\subset E(\c^m)$ denote the orthogonal complement: this is a vector
bundle of rank $m-1$ on which the quadratic form restricts with rank $m-2$, and
the quotient $Q = F^{\perp}/F$ is thus an $SO_{m-2}$-bundle, fitting into an
exact diagram:
$$
\begin{array}{rcccccccl}
&&0&&0&&&&\\
&&\downarrow&&\downarrow&&&&\\
&&&&&&&&\\
&&F&=&F&&&&\\
&&&&&&&&\\
&&\downarrow&&\downarrow&&&&\\
&&&&&&&&\\
0&\rightarrow&F^{\perp}&\rightarrow&E(\c^m)&\rightarrow&F^{-1}&\rightarrow&0\\
&&&&&&&&\\
&&\downarrow&&\downarrow&&\|&&\\
&&&&&&&&\\
0&\rightarrow&Q&\rightarrow&(F^{\perp})^{\vee}&\rightarrow&F^{-1}&\rightarrow&0
\\
&&&&&&&&\\
&&\downarrow&&\downarrow&&&&\\
&&0&&0&&&&\\
\end{array}
$$
Noting now that $E$ is determined up to isomorphism by $F^{\perp}$ together
with its (degenerate) quadratic form (which determines $F$), we see that to
construct such a bundle $E$ it is enough to specify the left-hand vertical
sequence. For this, $F$ is determined by $g$ parameters; $Q$ by
$(g-1)(m-2)(m-3)/2$ parameters; and the extension $F^{\perp}$ by
$$
h^1(C,Q\dual \otimes F) -1 = (m-2)(g-1-d) -1
$$
parameters. Consequently, a {\it generic} $E\in \mm(SO_m)$ (in either component
of the moduli space) possesses an isotropic line subbundle $F$ of degree $d$
only if the total number of parameters is at least $\dim \mm(SO_m) =
(g-1)m(m-1)/2$:
$$
\begin{array}{rcl}
g+ (g-1)(m-2)(m-3)/2\ \ \ \ \ \ \ &&\\
+ (m-2)(g-1-d) -1 &\geq& (g-1)m(m-1)/2,\\
\end{array}
$$
which simplifies to $d\leq -g +1$. The lemma follows at once from this.
\qed

\begin{rem}
\label{rank2case}\rm
Note that in the case $m=3$, $E(\c^3) = \ad V$, the bundle of tracefree
endomorphisms of a stable rank 2 vector bundle $V$---i.e. that coming from $E$
via the 2-dimensional spin representation. Then there is a one-to-one
correspondence between line subbundles $L\subset V$ and isotropic line
subbundles $F\subset \ad V$: this is because the quadratic form on $\ad V$ is
the Killing form, which for tracefree $2\times 2$ matrices is the determinant.
Thus $F$ is isotropic if and only if it consists of nilpotent endomorphisms,
and $L$ is then the kernel bundle; and conversely $F = {\rm Hom}(V/L, L) = L^2
\otimes \det V^{\vee}$.

If one views $L\subset V$ as determining a cross-section $ l \subset \bp (V)$
of the corresponding ruled surface, then $-\deg F = l\cdot  l$ is its
self-intersection. So lemma \ref{nagata} says
$$
l\cdot l \geq g-1
\qquad
\hbox{for any section $ l \subset \bp (V)$,}
$$
for a generic rank 2 vector bundle $V$. A well-known result of Nagata \cite{N},
on the other hand, says that {\it every} ruled surface has a section $l$ with
self-intersection $ l\cdot  l \leq g$. Since, for a given surface,
self-intersection of a section is constant ($\equiv \deg V$) mod 2, this
implies that the inequality of lemma \ref{nagata} is sharp for $m=3$.

\end{rem}

\subsection{Conservation of parity}

We shall need to make use of the following well-known result of Atiyah \cite{A}
and Mumford
\cite{M}:

\begin{mumlemm}
\label{mumlemm}
Suppose that $F \rightarrow C$ is a vector bundle admitting a nondegenerate
symmetric bilinear form
$$
F\otimes F \rightarrow K_C.
$$
Then $h^0(C,F)$ is constant modulo 2 under deformation.
\end{mumlemm}

The case of this which we shall be interested in arises when $F = L\otimes
E(\c^m)$ for $L$ a theta characteristic, $L^2 =K$, and $E$ an $SO_m$- (or
$SC_m$-) bundle. Let $\vartheta(C) \subset J^{g-1}(C)$ denote the set of theta
characteristics. The value of $h^0(C,F)$ mod 2 is in this case given by a
calculation of Serre:

\begin{prop}
\label{2.11}
For any $SC_m$-bundle $E\rightarrow C$ and theta characteristic $\th\in
\vartheta(C)$
one has
$$
h^0(C,\th \otimes E(\c^m)) \equiv m h^0(C,\th) + \deg \nm(E)
\quad {\rm mod}\ 2.
$$
\end{prop}

{\it Proof.}  By \cite{S}, theorem 2, one has the congruence:
$$
h^0(C,\th \otimes E(\c^m)) \equiv (m+1) h^0(C,\th) + h^0(C,\th \otimes w_1(E))+
w_2(E)
\quad {\rm mod}\ 2,
$$
where $w_1$ and $w_2$ are the Stiefel-Whitney classes. But $w_1(E)$ can be
identified with $\det E(\c^m) \in J_2(C) \cong H^1(C,\z/2)$, which in our case
vanishes since $E$ is a {\it special} Clifford bundle. On the other hand,
$w_2(E)\equiv \deg \nm(E)$ mod 2 by remark \ref{w2}. So we get the statement in
the proposition.
\qed

\begin{cor}
\label{serre}
Suppose that $p: \cctil \rightarrow C$ is an unramified double cover, and let
$\sigma: \cctil\leftrightarrow \cctil$ denote the sheet-interchange over $C$.
Suppose that
 $E$ is any $SC_m$-bundle on $C$, and that $L\in \pic(\cctil)$ satisfies
$L\otimes \sigma(L) = K_{\cctil}$. Then:
$$
h^0(\cctil,L \otimes p^* E(\c^m)) \equiv m h^0(\cctil,L)
\quad {\rm mod}\ 2.
$$
\end{cor}

\pf $L$ has degree $g(\cctil) -1$, so by the Atiyah-Mumford lemma it suffices
to assume that $L$ is the pull-back of a theta characteristic. Then the
corollary follows at once from the proposition since $w_2(p^* E)=0$.
\qed

\subsection{Generically zero results}




The main aim of this section is to establish that in each of the results
\ref{2.11} and \ref{serre} above, the vanishing (mod 2) of the right-hand side
guarantees that $h^0(C,L\otimes E(\c^m)) =0$ for generic stable orthogonal
bundle $E$ (theorems \ref{gen0}, \ref{gen0th} and corollary \ref{serregen0}).
To this end we shall
apply Brill-Noether methods (see, for example, \cite{ACGH}) to the loci:
$$
\ww_k = \{ (L,E) | h^0 (C, L\otimes E(\c^m)) \geq k \} \subset
J^{g-1}(C) \times \nn^{st}(m)
$$
where $\nn^{st}(m) \subset \nn(m)$ is the open set of {\it stable} bundles.
Analogously with classical Brill-Noether theory we claim that the Zariski
tangent space to $\ww_k$ at a point $(L,E)$ is annihilated by the image of a
{\it Petri map}
$$
\begin{array}{rcl}
\mu : H^0(L\otimes E(\c^m)) \otimes H^0(KL^{-1}\otimes E(\c^m))
&\rightarrow&
H^0(K) \oplus H^0(K\otimes E(\som))\\
&&\\
&&\cong \Omega^1_{J^{g-1}\times \nn^{st}(m)}|_{(L,E)}.\\
\end{array}
$$

Here $\som$ is the Lie algebra of $\spin_m$, viewed as the semisimple component
of the adjoint representation of $SC_m$. So $E(\som)$ is then the vector bundle
with fibre $\som$ associated to any Clifford bundle $E$ via this
representation; and by standard deformation theory (and this is where we
require our bundles to be stable)
$H^0(C,K\otimes E(\som)) \cong \Omega^1_{\nn^{st}(m)}|_{E}$.

The Petri map is defined by
\begin{equation}
\label{mudef}
\mu : s\otimes t \mapsto \< s,t\> \oplus s\wedge t
\end{equation}
where $\<,\>$ is the symmetric bilinear form on the vector bundle $E(\c^m)$,
and we identify $\som \cong \bigwedge^2 \c^m$.

\begin{prop}
\label{petri}
For $(L,E) \in \ww_k - \ww_{k+1}$ the Zariski tangent space of $\ww_k$ at
$(L,E)$ is
$
T_{(L,E)}\ww_k = (\im \mu)^{\perp}.
$
\end{prop}

To prove this we first need:

\begin{lemm}
Given $(L,E)\in \pic(C)\times \nn^{st}(m)$, a section $s\in
H^0(C,L\otimes E(\c^m))$, and a tangent vector $\eta \oplus \xi \in
H^1(C,\oo) \oplus H^1(C,E(\som)) = T_{(L,E)}\bigl( J^{g-1}\times \nn^{st}(m)
\bigr)$,
$s$ extends to the 1st-order deformation of $L \otimes E(\c^m)$
corresponding to $\eta \oplus \xi$ if and only if
$$
\xi  s +\eta s=0\in
H^1(C,L\otimes E(\c^m)).
$$
\end{lemm}

\pf
Represent $E$ by transition data $\{h_{\alpha \beta}\}$
with respect to an open cover $\{U_{\alpha}\}$ of $C$, where the $h_{\alpha
\beta}$
are holomorphic $SC_m$-valued functions on $U_{\alpha}\cap
U_{\beta}$ satisfying the cocycle condition; let $\gab$ be the image of
$h_{\alpha \beta}$ in $SO_m$---these are then the transition functions for the
vector bundle $E(\c^m)$.
Likewise represent $L$ by
transition data $\{ \phab \}$ where the $\phab$ are $\c^*$-valued
functions. Finally,
represent the sum
$\eta \oplus \xi \in H^1(\oo) \oplus
H^1(\som)$ by a cocycle $\{\eab \oplus \xab\}$ (holomorphic Lie
algebra-valued functions on the $U_{\alpha}\cap U_{\beta}$).

Then the
1st-order deformation corresponding to $\eta \oplus \xi$
is the $SO_m$-bundle on
$C\times \spec \c[\varepsilon ] /(\varepsilon^2)$ with transition data
$$
\begin{array}{rcl}
\tilde \phab &=&  (1+\varepsilon \eab)\phab,\\
\tilde \gab &=&  (1+\varepsilon \xab)\gab.\\
\end{array}
$$
A section $s\in H^0(C,L\otimes E(\c^m)) $ is now given by a collection
$\{s_{\alpha}\}$ of
$\c^m$-valued functions satisfying
$$
s_{\alpha} = \phab \gab s_{\beta}
$$
on $U_{\alpha}\cap U_{\beta}$; and $s$ extends to the 1st-order
deformation $\xi$ provided there exists a 0-cochain $\{s_{\alpha}'\}$ such
that
$$
\tilde s_{\alpha} = s_{\alpha} + \varepsilon s'_{\alpha}
$$
satisfies the cocycle condition
$$
\begin{array}{rcl}
\tilde s_{\alpha} &=& \tilde \phab \tilde \gab  \tilde s_{\beta} \\
   &=& \phab (1+\varepsilon \eab)(1+\varepsilon \xab )\gab (s_{\beta} +
\varepsilon s'_{\beta})\\
    &=& s_{\alpha} + \varepsilon  (\eab s_{\alpha} +\xab s_{\alpha}
   +\phab \gab s'_{\beta} ),\\
\end{array}
$$
and hence
$
s'_{\alpha} - \phab \gab s'_{\beta} = (\eab + \xab)s_{\alpha}.
$
In other words $s'= \{s'_{\alpha}\}$ has to satisfy
$$
\eta s + \xi s = \partial s' \in Z^1(L \otimes E(\c^m)),
$$
where $\partial$ is the coboundary operator on 0-cocycles, and $Z^1$ is
the group of Cech 1-cocycles. So there exists a solution if and only if
$\xi s +\eta s=0$ in cohomology.
\qed

{\it Proof of proposition \ref{petri}.}
We write down the Serre duality pairing of $\eta \oplus \xi \in  H^1(C,\oo)
\oplus H^1(C,E(\som))$ with $\mu (s\otimes t)$, for $s\in H^0(C, L\otimes
E(\c^m))$ and $t\in H^0(C, KL^{-1}\otimes E(\c^m))$:
\begin{equation}
\label{pair}
\begin{array}{rcl}
\bigl(\eta \oplus \xi,\mu (s\otimes t)\bigr)_{\rm Serre}
&=&\displaystyle
\int_C \bigl(\eta \<s,t\> + {\rm trace}(\xi\ s\wedge t)\bigr)\\
&&\\
&=&\displaystyle \int_C \< \xi s+\eta s,t\>,\\
\end{array}
\end{equation}
where we have identified $\bigwedge^2 \c^m$ with $\som$, the space of
skew-symmetric $m\times m$ matrices, and used the fact that under this
identification one has ${\rm trace}(\xi\ s\wedge t) = \<\xi s,t\>$, as one
verifies by an easy calculation.

Now for $(L,E) \in \ww_k - \ww_{k+1}$, the Zariski tangent space
$T_{(L,E)}\ww_k$
is the linear span of directions $\eta \oplus \xi$ in which {\it all} sections
$s$ extend. By the lemma this condition on $\eta \oplus \xi$ is that $\xi s +
\eta s =0$ in cohomology for all $s\in H^0(C, L \otimes E(\c^m))$. By
(\ref{pair}) this linear span is precisely $(\im \mu)^{\perp}$.
\qed

We come now to the first main result of this section:

\begin{theo}
\label{gen0}
For generic $(L,E)\in J^{g-1}(C)\times \mm^{\pm}(SO_m)$ one has $H^0(C,L\otimes
E(\c^m)) =0$.
\end{theo}

\pf
We choose a component of $J^{g-1} \times \nn^{st}(m)$, and suppose that
$h^0(C,L\otimes E(\c^m))\geq k>0$ everywhere in this component, and generically
equal to $k$, for some natural number $k$. We then choose a generic point
$(L,E)$ in this component; according to proposition \ref{petri} the Petri map
$\mu$ is identically zero here. This means that for arbitrary sections $s\in
H^0(C,L\otimes E(\c^m))$ and $t\in H^0(C,KL^{-1}\otimes E(\c^m))$ (and note
that both spaces have the same dimension $k$, by Riemann-Roch) we have, on the
one hand, $s\wedge t =0$. This implies that $s,t$ generically generate the same
line subbundle $F\subset E(\c^m)$; moreover, since $s,t$ are arbitrary we see
that $F$ is independent of their choice, and depends only on $L$ and $E$.

On the other hand, again by (\ref{mudef}), $\<s,t \>=0$; this implies that
$F$ is an {\it isotropic} line subbundle. (Of course, since $E$ is generic this
is also forced by corollary \ref{alliso}.)
Lemma \ref{nagata} and genericity of $E$ therefore implies that
$
\deg F \leq -g+1.
$
But by definition of $F$ the spaces of sections $H^0(L\otimes F)$ and
$H^0(KL^{-1}\otimes F)$ are both nonzero, which
forces $L\otimes F = KL^{-1}\otimes F = \oo$ and hence
$L^2 = K$. But this contradicts the genericity of $L$.
\qed

We now consider the intersection of the subscheme $\ww_k$ with the fibre
$\{\th\} \times \nn^{st}(m)$, where $\th \in \vartheta(C)$ is a theta
characteristic. (Indeed, theorem \ref{gen0} above also follows---for all cases
except $\mm^-(SO_{2n})$---from theorem \ref{gen0th} below.)

Fixing $\th$,
one sees that by Serre duality the Petri map now factorises into two maps:
$$
\begin{array}{rcl}
\mu_J : S^2 H^0(C,\th \otimes E(\c^m)) & \rightarrow & H^0(C,K),\\
\mu_{\nn}: \bigwedge^2 H^0(C,\th \otimes E(\c^m)) & \rightarrow &
H^0(C,K\otimes\bigwedge^2  E(\c^m))\\
&& \cong H^0(C,K\otimes E(\som)),\\
\end{array}
$$
where $\mu_J : s\otimes t \mapsto \<s,t\>$, and $\mu_{\nn}$ is the natural
multiplication map.

As an immediate consequence of proposition \ref{petri} we obtain:

\begin{cor}
\label{thpetri}
If $E\in \nn^{st}(m)$ and $\th \in \vartheta(C)$ satisfy $h^0(C,\th \otimes
E(\c^m)) = k$ then the subscheme $\uu_k = \ww_k |_{\{K^{1/2}\} \times
\nn^{st}(m)} \subset \nn^{st}(m)$ has Zariski tangent space
$
T_E \uu_k = (\im \mu_{\nn})^{\perp}.
$
\end{cor}

{}From this follows our second main result:

\begin{theo}
\label{gen0th}
For $m\geq 3$ and for
any theta characteristic $\th \in \vartheta(C)$ we have $h^0(C,\th \otimes
E(\c^m)) = 0$ or 1 for generic $E\in \mm^{\pm}(SO_m)$, where the parity is
determined by lemma \ref{2.11}.
\end{theo}

\pf
We fix our component $\mm^{\pm}(SO_m)$, and suppose that $h^0(C,\th\otimes
E(\c^m))\geq k$ everywhere in this component, generically equal to $k$. We then
choose a generic point $E$ in this component; by proposition \ref{thpetri} the
Petri map $\mu_{\nn}$ vanishes here. As in the proof of theorem \ref{gen0} this
means $s\wedge t =0$ for all sections $s,t \in H^0(C,\th\otimes E(\c^m))$, so
all sections generate the same line subbundle $F\subset E(\c^m)$; in particular
$$
H^0(C,\th\otimes F) =
H^0(C,\th\otimes E(\c^m)).
$$
Since $E$ is generic it follows by corollary \ref{alliso} that $F$ is
isotropic,
and hence by lemma \ref{nagata} that $\deg F \leq -g+1$. But then $\deg
\th\otimes F \leq 0$ and hence $k=h^0(C,\th\otimes F) \leq 1$.
\qed

\begin{cor}
\label{serregen0}
Suppose that $p: \cctil\rightarrow C$ is an unramified double cover and that
$E,L$ is a {\it generic} pair as in corollary \ref{serre}. Then---except
possibly in the case of $m$ even, $w_2(E) \equiv 1$---we have
$
h^0(\cctil,L \otimes p^* E(\c^m)) =0
$ or 1, where the parity is determined by corollary \ref{serre}.
\end{cor}

\pf
As in the proof of \ref{serre} we take $L= p^* \th$ to be the pull-back of a
theta characteristic. (See also the discussion of section \ref{main}.) Then
$$
h^0(\cctil,L \otimes p^* E(\c^m)) = h^0(C,\th\otimes E(\c^m)) +
h^0(C,\th\otimes\eta \otimes E(\c^m)),
$$
where $\eta\in J_2(C)$ is the 2-torsion point associated to the covering.
If $E$ is generic then by theorem \ref{gen0th} it is possible to choose
$\th$---in all cases except when $m$ is even and $w_2(E)$ is odd---so that the
right-hand side is $\leq 1$.
\qed

\begin{rem}\rm
In the case of $m$ even, $w_2(E)$ odd, $h^0(\cctil,L \otimes p^* E(\c^m))\geq
2$ for any $L=p^*\th$; so the above argument fails. Nonetheless, we expect that
the result \ref{serregen0} is still true in this case, but its proof requires a
refinement of the Brill-Noether analysis of this section.
\end{rem}

\section{The Verlinde formula}
\label{vformula}

In this section we shall write down, for the unitary and spin groups, the
Verlinde formula which calculates the dimension of the vector spaces
$H^0(\mm(G),\Theta(V))$.
For the derivations of these formulae we refer the reader to \cite{B3},
\cite{OW}.
In fact, what one writes down is a natural number $N_l(G)$
depending on the group, on the genus $g$, and on an integer $l$ called the
`level'. Then to any representation $V$ of $G$, one associates a level $l=d_V$
such that
$\dim H^0(\mm(G),\Theta(V)) = N_{d_V}(G)$.

\subsection{Preliminaries}

For $k\in \n$ and $r\in \q$ we let
$$
f_k(r)
=4\sin^2(r\pi/k)=(1-\zeta_{k}^r)(1-\zeta_{k}^{-r})
$$
where $\zeta_k = e^{2\pi i /k}$. This satisfies certain obvious identities (we
shall usually drop the subscript $k$ for convenience):

\begin{lemm}
\label{fids1}
\item[(i)] $f(r)=f(-r)$;
\item[(ii)]$ f(r)=f(k-r)$;
\item[(iii)] $ f(k/2)=4$;
\item[(iv)] $ f(2r)=f(r)f(k/2-r)$;
\item[(v)] $ \prod_{r=1}^{k-1}f(r)=k^2$.
\end{lemm}

In addition, we shall need the following:

\begin{lemm}
\label{fids2}
\item[(i)]
If $k=2n+1$ is odd then
$$
\begin{array}{lcl}
 \prod_{r=1}^{n}f(r) &=& 2n+1,\\
 \prod_{r=1}^{2n}f(r)^{[{r\over 2}]} &=& (2n+1)^n.\\
\end{array}
$$
\item[(ii)]
If $k=2n$ is even then
$$
\begin{array}{lcl}
 \prod_{r=1}^{n-1}f(r) &=& n,\\
 \prod_{r=1}^{2n-1}f(r)^{[{r\over 2}]} &=& 2^{n-1} n^n\\
 \prod_{r=1}^{2n-1}f(r)^{[{r+1\over 2}]} &=& 2^{n+1} n^n.\\
\end{array}
$$
\end{lemm}

\pf (i) The first identity follows at once from parts (ii) and (v) of the
previous lemma. For the second, use \ref{fids1} (ii) to write
$$
\prod_{r=1}^{2n}f(r)^{[{r\over 2}]} = \prod_{r=1}^{n}f(r)^{[{r\over 2}] +
[{2n+1-r\over 2}]} = \prod_{r=1}^{n}f(r)^n = (2n+1)^n.
$$

(ii) Again the first identity is an easy consequence of the previous lemma.
We shall just give the proof of the third identity, that of the second being
being almost the same.

The left-hand product can be rewritten, using \ref{fids1} (ii), as:
$$
P =
f(n)^{[{n+1\over 2}]} \times
\prod_{r=1}^{n-1}f(r)^{[{r+1\over 2}] + [{2n-r+1\over 2}]}.
$$
We observe that
$$
\Bigl[{r+1\over 2}\Bigr] + \Bigl[{2n-r+1\over 2}\Bigr] =
\cases{n & if $r$ is even,\cr
       n+1 & if $r$ is odd.\cr}
$$
So suppose first that $n$ is even. Using \ref{fids1} (iii) and
$\prod_{r=1}^{n-1}f(r) = n$ we see that
$$
\begin{array}{rcl}
P &=& 2^n f(1)^{n+1}f(2)^n \cdots f(n-2)^n f(n-1)^{n+1}\\
 &=& (2n)^n f(1)f(3) \cdots f(n-1).\\
\end{array}
$$
We claim that $f(1)f(3) \cdots f(n-1) = 2$---from which the third identity in
\ref{fids2} (ii) follows. To see this, observe first that (using \ref{fids1}
(iv))
$$
\begin{array}{rcl}
f(2)f(4) \cdots f(n-2) &=& f(n/2)^{-1} \prod_{r=1}^{n-1}f(r) \\
                       &=& n/2;\\
\end{array}
$$
and second that
$$
n = f(1) \cdots f(n-1) \times  f(2) \cdots f(n-2).
$$
The reasoning for $n$ odd is similar.
\qed

In what follows
we shall consider sets $U = \{u_1,\ldots ,u_n\}$ of rational numbers; and for
such a set we define:
$$
\begin{array}{rcl}
\Psi_k(U) &=& \displaystyle\prod_{1\le i<j\le n}  f(u_i -u_j) \\
&&\\
\Pi_k(U) &=& \displaystyle\prod_{1\le i<j\le n}  f(u_i -u_j)f(u_i +u_j) \\
\Phi_k(U) &=& \displaystyle \prod_{1\le i<j\le n}  f(u_i -u_j)f(u_i +u_j)
\times
\prod_{i=1}^n f(u_i)
.\\
\end{array}
$$

\subsection{Computations}

\begin{ex}
\label{4.1} $\bf SL_{n}.$ \rm
The Verlinde number in this case is:
$$
N_l(SL_n) = \sum_U \Bigl( {n(l+n)^{n-1}\over \Psi_{l+n}(U)}
\Bigr)^{g-1}
$$
where the sum is taken over $U = \{ 0=u_0<u_1<\cdots <u_{n-1}<l+n \}$.

Since
we are concerned with the spin groups we shall record at this point the result
of computing this expression at levels 1 and 2 for $\spin_6 = SL_4$:
\begin{equation}
\label{spin6}
\begin{array}{crcl}
\bullet &   N_1(SL_4) &=& 2^{2g};\\
\bullet &    N_2(SL_4) &=& 2^{3g-1}3^{g-1} + 2^{3g-1} + 2^g 3^{g-1}.\\
\end{array}
\end{equation}

\end{ex}

\begin{ex}
\label{4.2}
 $\bf Spin_{2n},\ n \geq 4.$ \rm
The Verlinde number at level $l$ is here:
$$
N_l(\spin_{2n}) = \sum_{U\in P_l(2n)} \Bigl( {4k^n \over \Pi_k (U)}
\Bigr)^{g-1},
\qquad {\rm where}\ k=l+2n-2;
$$
where $P_l = P_l(2n)$ denotes the collection of sets $U=\{u_1<\cdots <u_n\}$
satisfying the following conditions. (Note that this collection is finite
provided $n\geq 3$.)
\begin{enumerate}
\item  $ u_i \in {1\over 2}\z$;
\item  $ u_{i+1} -u_{i} \in \z$ for $i=1,\ldots n-1$;
\item  $ u_1 +u_2 >0$; and
\item $ u_{n-1} +u_n < k=l+2n-2$.
\end{enumerate}

It will be convenient to write $P_l = P_l^+ \cup P_l^-$ where $P_l^+$ (resp.
$P_l^-$) consists of $U$ with all $u_i \in \z$ (resp. $u_i \in {1\over
2}\z\backslash \z$). Correspondingly the Verlinde number splits up as $N_l =
N_l^+ + N_l^-$ where
$$
N_l^{\pm}(\spin_{2n}) = \sum_{U\in P_l^{\pm}(2n)} \Bigl( {4k^n \over \Pi_k (U)}
\Bigr)^{g-1}.
$$

\begin{rem} \rm
$P_l(2n)$ may be viewed as a set of highest weights of irreducible
representations of $\spin_{2n}$ (namely, those weights in a fundamental chamber
for the action of the affine Weyl group of level $l$). Then $P_l^+$ is the
subset of `tensor' representations---those which descend to $SO_{2n}$---and
$P_l^-$ is the subset of `spinor' representations. The same remark applies to
the odd spin groups below.
\end{rem}

As for $SL_4$ let us note the lowest cases of this formula.
(Note, incidentally, that the formulae (\ref{spin2n}) are consistent with
$\spin_4 = SL_2\times SL_2$ and $\spin_6 = SL_4$: for $n=2,3$ respectively
$N_l$ coincides with $N_l(SL_2)^2$ and $N_l(SL_4)$.)
\begin{equation}
\label{spin2n}
\begin{array}{clcl}
\bullet &  N_1(\spin_{2n}) &=& 2^{2g}\\
\bullet &  N_2^+(\spin_{2n}) &=& (2n)^{g-1}(n-1) + 2^{3g-1}n^{g-1} \\
\bullet &  N_2^-(\spin_{2n}) &=& 2^{3g-1}.\\
\end{array}
\end{equation}

\pf
We shall just prove the formulae for $N_2^{\pm}$, and leave $N_1$ (which we
shall not need) to the reader.

We begin by listing the sets over which the summation takes place:
$$
\begin{array}{lcl}
P_2^+(2n): && \{0,1,\ldots, n-2 ,n\pm 1 \} \\
           && \{0,1,\ldots, n-2,n \} \\
           && \cdots \\
           && \{0,2,\ldots, n-1,n \} \\
           && \{\pm 1,2,\ldots, n-1,n \}; \\
&&\\
P_2^-(2n): && \{\pm{1\over 2},{3\over 2}, \ldots, n \pm {1\over 2}  \}. \\
\end{array}
$$
We recall from \cite{OW} that reflection of the end-points in $0,n$---i.e. the
$\pm$ signs in these sets---defines an action of $\z/2 \times \z/2$ on
$P_2(2n)$ under which $\Pi_k(U)$ is easily seen to be invariant. Thus, for
example, $P_2^-(2n)$ is a single orbit and so
\begin{equation}
\label{N2-}
N_2^-(\spin_{2n}) = 4 \times \Bigl({4(2n)^n\over \Pi_{2n}(U)}\Bigr)^{g-1}
\end{equation}
where
$
U = \{{1\over 2},{3\over 2} \ldots,n - {1\over 2}  \}.
$
We have here
$$
\begin{array}{rcl}
\Pi_{2n}(U) &=& \displaystyle
\prod_{1\leq i<j\leq n } f_{2n}(i-j) f_{2n}(i+j-1) \\
            &=& \displaystyle\prod_{r=1}^{2n-1} f_{2n}(r)^{m_r},\\
\end{array}
$$
where the multiplicities $m_r$ are to be determined.
Namely, $m_r = a_r + b_r$ where $a_r$ is the number of pairs $i<j$ such that
$j-i = r$, and $b_r$ is the number of pairs $i<j$ such that $i+j = r+1$. From
this we see that
$$
\begin{array}{rcl}
a_r &=&
\cases{n-r & if $r<n$\cr
       0 & if $r\geq n$,\cr}\\
b_r &=&
\cases{[ {r/ 2}] & if $r<n$\cr
       [ {(2n-r)/ 2}] & if $r\geq n$;\cr}\\
\end{array}
$$
and hence that
$$
m_r = n-\Bigl[ {r+1\over 2}\Bigr].
$$
So it follows from \ref{fids1} (v) and \ref{fids2} (ii) that
$$
\Pi_{2n}(U) = {(2n)^{2n}\over \prod_{r=1}^{2n-1}f(r)^{[{r+1\over 2}]}}
={1\over 2}(2n)^n ;
$$
and hence from (\ref{N2-}) that $N^-_2(\spin_{2n}) = 2^{3g-1}$.

\medskip

Let us now turn to $N^+_2(\spin_{2n})$. For $l=0,1,\ldots ,n$ we shall write
$$
U_l = \{ 0,1,\ldots,\widehat l,\ldots ,n\},
$$
i.e. $l$ has been deleted. Thus $P_2^+(2n)$ consists of $U_0,U_1,\ldots,U_n$
together with the reflections of $U_0$ and $U_n$ under the action of $\z/2
\times \z /2$; and so
\begin{equation}
\label{N2+}
\begin{array}{rcl}
N_2^+(\spin_{2n}) &=& \displaystyle
2 \Bigl( {4(2n)^n \over \Pi_{2n}(U_0)} \Bigr)^{g-1}
+2 \Bigl( {4(2n)^n \over \Pi_{2n}(U_n)} \Bigr)^{g-1}\\
&&\\
&&\displaystyle\ \ \ \ \ \ \
+\sum_{l=1}^{n-1}\Bigl( {4(2n)^n \over \Pi_{2n}(U_l)} \Bigr)^{g-1}.\\
\end{array}
\end{equation}
Now we can write, for each $l$, $\Pi_{2n}(U_l) = N/D_l$ where:
$$
\begin{array}{rcl}
N &=& \displaystyle
\prod_{0\leq i<j\leq n} f(j-i)f(i+j),\\
D_l &=&\displaystyle \prod_{j=l+1}^n f(j-l)f(j+l) \prod_{i=0}^{l-1}
f(l-i)f(l+i)
\quad {\rm for}\ 1\leq l \leq n-1.\\
\end{array}
$$
For $l=0$ and $n$ the denominator takes the slightly simpler forms:
$$
D_0 = \prod_{j=1}^n f(j)^2 = 16 n^2,
$$
using \ref{fids1} (iii) and \ref{fids2} (ii); and likewise
$$
\begin{array}{rcl}
D_n &=& \prod_{i=0}^{n-1} f(n-i)f(n+i)\\
    &=&\prod_{i=0}^{n-1} f(n-i)^2\\
 &=&\prod_{i=1}^{n} f(i)^2\\
 &=& 16 n^2.\\
\end{array}
$$
To compute the denominator $D_l$ for $1\leq l \leq n-1$:
$$
\begin{array}{rclr}
D_l &=& f(1) \cdots f(n-l) \times f(2l+1)\cdots f(n+l)&\\
    && \times f(1) \cdots f(l) \times f(l) \cdots f(2l -1)&\\
&&&\\
&=&\displaystyle f(n){f(l)\over f(2l)} \prod_{r=1}^{n-1}f(r)^2 {f(n+1)\cdots
f(n+l)\over f(n-1)\cdots f(n-l+1)}&\\
&&&\\
&=& 4n^2 {f(l)f(n+l)/ f(2l)}&\hbox{by \ref{fids1} (ii),(iii)}\\
&&&\hbox{and \ref{fids2} (ii);}\\
&=& 4n^2 {f(n+l)/ f(n-l)}&\hbox{using \ref{fids1} (iv);}\\
&=& 4n^2.&\\
\end{array}
$$
We next compute the numerator:
$$
N = \prod_{r=1}^{2n-1} f(r)^{m_r}
$$
where $m_r = a_r + b_r$ and $a_r$ is the number of pairs $i<j$ between $0$ and
$n$ such that $j-i =r$, and $b_r$ the number such that $i+j =r$. So:
$$
\begin{array}{rcl}
a_r &=&
\cases{n+1-r & if $r\leq n$\cr
       0 & if $r> n$,\cr}\\
b_r &=&
\cases{[ {(r+1)/ 2}] & if $r\leq n$\cr
       [ {(2n-r+1)/ 2}] & if $r> n$;\cr}\\
\end{array}
$$
from which we find
$$
\begin{array}{rcl}
m_r &=&
\cases{n+1 -[r/2] & if $r\leq n$\cr
       n -[r/2] & if $r> n$.\cr}\\
\end{array}
$$
Hence
$$
N = \prod_{r=1}^{2n-1} f(r)^{n-[r/2]} \times \prod_{r=1}^n f(r) = 4(2n)^{n+1},
$$
using \ref{fids1} (v) and \ref{fids2} (ii).

Putting together the above computations we obtain
$
\Pi_{2n}(U_l) = (2n)^{n-1}$ for $l=0,n$, and $2^{n+1}n^{n-1}$ for $l=1,\ldots
,n-1$. Substituting into (\ref{N2+}) this gives
$$
N^+_2(\spin_{2n}) = (2n)^{g-1}(n-1) + 2^{3g-1}n^{g-1}
$$
as asserted.
\qed

\end{ex}

\begin{ex}
\label{4.3}
$\bf Spin_{2n+1}, n \geq 2.$ \rm
In this case the Verlinde number is:
$$
N_l(\spin_{2n+1})= \sum_{U\in P_l(2n+1)} \Bigl( {4k^n \over \Phi_k(U)}
\Bigr)^{g-1},
\qquad {\rm where}\ k=l+2n-1;
$$
where $P_l(2n+1)$ consists of $U=\{0<u_1 <\cdots < u_n \}$ such that
\begin{enumerate}
\item  $ u_i \in {1\over 2}\z$;
\item  $ u_{i+1} -u_{i} \in \z$ for $i=1,\ldots n-1$;
\item $ u_{n-1} +u_n < k=l+2n-1$.
\end{enumerate}

As in the even case we shall write
$$
N_l^{\pm}(\spin_{2n+1}) = \sum_{U\in P_l^{\pm}} \Bigl( {4k^n \over \Phi_k (U)}
\Bigr)^{g-1},
$$
where $P_l^{\pm}$ denotes the subsets of integral and half-integral $U$
respectively.

Here the lowest Verlinde numbers are:
\begin{equation}
\label{spin2n+1}
\begin{array}{clcl}
\bullet &  N_1(\spin_{2n+1}) &=& 2^{g-1}(2^g + 1)\\
\bullet &  N_2^+(\spin_{2n+1}) &=& 2^{2g-1} \\
\bullet &  N_2^-(\spin_{2n+1}) &=& (2n+1)^{g-1}(2^{2g-1} + n).\\
\end{array}
\end{equation}

We shall omit the proof of these formulae, as it is entirely similar to that
of (\ref{spin2n}). Moreover they are proved (by a somewhat clumsier method) in
\cite{O}.
\end{ex}

\subsection{Sections of theta bundles}

Our interest in the preceding calculations lies in the fact that $N_2(\spin_m)$
is the dimension of $H^0(\mm(\spin_m),\Theta(\c^m))$, where $\Theta(\c^m)$ is
the theta line bundle, defined in section \ref{mod}, for the standard
orthogonal representation. To each irreducible representation $V$ of a group
$G$ there is associated an integer $d_V$---the {\it height} or {\it Dynkin
index} of the representation---for which, for any $r\in \z$, one has

\begin{equation}
\label{verlinde}
h^0(\mm(G),\Theta(V)^{\otimes r}) = N_{rd_V}(G).
\end{equation}
(For a useful discussion of the Dynkin index and the proof of (\ref{verlinde})
see \cite{LS}; see also
\cite{BL}, \cite{F}, \cite{KNR}.)

When $V=\c^m$ is the standard orthogonal representation the
height is:
\begin{equation}
\label{di2}
d_V=
\cases{
4 & if $m=3$,\cr
2 & if $m\geq 5$.\cr}
\end{equation}

What is the analogue of (\ref{verlinde}) for the `twisted' moduli space
$\mm^-(\spin_m)$ defined in section \ref{clif}? The following is a refinement
(which deals with even as well as odd $m$) of conjecture (5.2) in \cite{OW}.

\begin{conj}
\label{conj}
For $m\geq 3$, and any representation $SC_m \rightarrow SL(V)$ satisfying
(\ref{ss}),
$$
h^0(\mm^-(\spin_m),\Theta(V)) = (-1)^m (N_{d_V}^+ -N_{d_V}^-).
$$
\end{conj}

\begin{ex} $\bf m=3.$ \rm
By example \ref{4.1},
$$
N_l(\spin_3) = N_l(SL_2) = \sum_{j=1}^{l+1} \Bigl( {l+2 \over 2\sin^2{j\pi
\over l+2}}
\Bigr)^{g-1}.
$$
Notice that formally computing
\ref{4.3} with $n=1$ (strictly speaking the Verlinde number \ref{4.3} makes
sense only for $n\geq 2$) we find (taking $u_0 = 0$):
$$
N_l(\spin_{2n+1})|_{n=1} = \sum_{j=1}^{2l+1} \Bigl( {l+1 \over 4\sin^2{j\pi
\over 2l+2}}
\Bigr)^{g-1} = N_{2l}(SL_2).
$$
Now, by (\ref{verlinde}) and the fact that $\ll_0 = \Theta(\c^2)$ where
$d_{\c^2} =1$, we have
$$
h^0(\mm(\spin_3), \Theta(V)) = h^0(\su(2),\ll_0^{d_V}) =
N^+_{d_V}(SL_2)+N^-_{d_V}(SL_2).
$$
(See example \ref{clifm=3}.)

On the other hand, on
$\mm^-(\spin_3) = \su(2,1)$ we have, given a Clifford representation $V$ for
which $d_V$ is even, $\Theta (V) = \ll_0^{d_V} = \ll_1^l$ where $2l = d_V$.
Thaddeus's twisted Verlinde formula (\cite{T}, corollary (18)) tells us that
$$
\begin{array}{rcl}
h^0(\su(2,1),\ll_1^{l}) &=& \displaystyle \sum_{j=1}^{2l+1} (-1)^{j+1}\Bigl(
{l+1 \over 4\sin^2{j\pi \over 2l+2}}
\Bigr)^{g-1}\\
&&\\
 &=& -N^+_{2l}(SL_2)+N^-_{2l}(SL_2);\\
\end{array}
$$
or in other words
$
\mm^-(\spin_3),\Theta(V)) =
-N^+_{d_V}(SL_2)+N^-_{d_V}(SL_2).
$

\end{ex}

\begin{ex} $\bf m=6.$ \rm
By examples \ref{clifm=6} and \ref{4.1}
$$
h^0(\mm(\spin_6), \Theta(V)) = h^0(\su(4),\ll_0^{d_V}) =
N^+_{d_V}(SL_4)+N^-_{d_V}(SL_4),
$$
where
$$
h^0(\su(4),\ll_0^{l}) =
N_l(SL_4) = \sum_{U} \Bigl( {4(l+4)^3 \over \Psi_{l+4}(U)}
\Bigr)^{g-1},
$$
summed over $U=\{0=u_0<u_1<u_2<u_3\leq l+3 \}\subset \z$.

On the other hand, $\mm^-(\spin_6) = \su(4,2)$ and by the same arguments as in
the previous example, the conjecture is in this case equivalent to:
$$
\begin{array}{rcl}
h^0(\su(4,2),\ll_2^{l}) &=& N^+_{2l}(SL_4) - N^-_{2l}(SL_4)\\
&&\\
&=& \displaystyle
\sum_{U} (-1)^{u_2}\Bigl( {4(2l+4)^3 \over \Psi_{2l+4}(U)}
\Bigr)^{g-1},\\
\end{array}
$$
summed over $U=\{0=u_0<u_1<u_2<u_3\leq 2l+3 \}\subset \z$.
One can possibly verify this using \cite{BL} theorem 9.4.
\end{ex}

\section{Numerology}
\label{numerology}

Our principal aim in this section to make some sense of the Verlinde numbers
(\ref{spin6}), (\ref{spin2n}) and (\ref{spin2n+1})
computed in the previous section at level 2---those, that is, associated to the
orthogonal representation of the Clifford group. But we begin with a remark,
independent of the rest of the paper, which deals with the Verlinde number
(\ref{spin2n+1}) at level 1.

\subsection{A remark on the generator of $\pic\ \mm(\spin_{2n+1})$}
\label{pfaffian}

It is known from \cite{LS} that the Picard group of the moduli scheme
$\mm(\spin_{m})$ is infinite cyclic. When $m=2n+1$ is odd the
ample generator $\pp$ is constructed as a Pfaffian bundle, i.e. $2\pp =
\Theta(\c^{2n+1})$. (When $m=2n$ is even $\pp$ exists only as a Weil divisor
class.)
Then according to
(\ref{spin2n+1}) the space of sections of this line bundle has dimension

$$
h^0(\mm(\spin_{2n+1}), \pp) = N_1(\spin_{2n+1}) = 2^{g-1}(2^g +1).
$$
This formula is striking, first because it is independent of $n$ and second
because it is the number of even theta characteristics of the curve.
In fact it is easy to see how to construct a basis for the linear system
$|\pp|$, as follows.

Let $\vartheta^+(C)$ denote the set of even theta characteristics, and for each
$L \in \vartheta^+(C)$ consider the reduced divisor
$$
D_L = \{ E \in \mm(\spin_{2n+1}) | H^0 (C,L\otimes E(\c^{2n+1})) \not=0\}.
$$
By proposition \ref{2.11} the dimension of $H^0 (C,L\otimes E(\c^{2n+1}))$ is
always even, and so by theorem \ref{gen0th} $D_L$ is a proper subset of the
moduli space. Consequently it is divisor---this is shown in proposition
\ref{div} below---and in fact, since we take $D_L$ to be reduced, some multiple
$k D_L \in |\Theta(\c^{2n+1})| = |2 \pp|$. Since the dimension of $H^0
(C,L\otimes E(\c^{2n+1}))$ is always even we see that $k$ is at least, and
hence equal to, 2.
We have therefore constructed a set of divisors on $\mm(\spin_{2n+1})$:
$$
D_L \in | \pp |,
\qquad
L\in \vartheta^+(C).
$$
In the case $n=1$ these were shown by Beauville \cite{B2} to be linearly
independent, and hence a basis of the linear system; we expect the same to be
true in general.

\subsection{Prym varieties}
\label{pryms}

We recall the following `Verlinde numbers' for principally polarised abelian
varieties. Let $(A,\Xi)$ be any principally polarised abelian variety of
dimension $g$, where $\Xi$ is a symmetric divisor representing the
polarisation;
and let
$$
H^0(A, m\Xi) = H^0_+(A, m\Xi)\oplus H^0_-(A, m\Xi)
$$
be the decomposition into $\pm$-eigenspaces under the canonical involution of
$A$. Then by writing down a suitable basis of theta functions one can easily
verify that:

\begin{equation}
\label{abtheta}
\dim H^0_{\pm}(A, m\Xi) =
\cases{
{(m^g \pm 2^g )/ 2} & if $m\equiv 0$ mod 2,\cr
{(m^g \pm 1 )/ 2} & if $m\equiv 1$ mod 2.\cr}
\end{equation}

Associated to a smooth projective curve $C$ we have a natural configuration of
principally polarised {\it Prym varieties}. Let us recall the usual notation
(see \cite{ACGH}).
 For each nonzero half-period $\eta \in J_2(C) \backslash
\{\oo\}$ we have an unramified double cover
$$p : \ctil \rightarrow C.$$
Writing $\jtil = J^{2g-2}(\ctil)$ we have
$$\nmp ^{-1}(K_C) = \peta \cup \peta^- \subset \jtil;$$
where $\peta,\ \peta^-$ are disjoint translates of the same abelian subvariety,
characterised by the condition that for $L\in \nmp^{-1}(K_C)$:

\begin{equation}
\label{prym+-}
h^0(\ctil, L) \equiv \cases{0& mod 2 if $L\in \peta$,\cr
                              1& mod 2 if $L\in \peta^-$.\cr}
\end{equation}
Then $\peta$ is called the Prym variety of the covering.
We shall denote by $\xieta$ the symmetric divisor representing the canonical
principal polarisation on $\peta$,
defined by $2 \xieta = \peta \cap \widetilde \Theta$, where $\widetilde \Theta$
is the theta-divisor in $\jtil$.
We shall allow also $\eta = 0$ by setting $(P_0 ,\Xi_0) = (J^{g-1}(C),\theta)$.

Notice that on $\peta^-$ both the polarisation and the line bundle $2\xieta$
are defined, and identify under translation with the same objects on $\peta$.
The distinguished line bundle $\xieta$ which represents the polarisation is no
longer defined on $\peta^-$, however. Thus in what follows the line bundle
$m\xieta$ makes sense on $\peta ^-$ {\it only if $m$ is even.}

\subsection{Spin Verlinde numbers versus Prym Verlinde numbers}

We may now formulate the computations of section \ref{vformula} as follows.

\begin{theo}
\label{numer}
\begin{enumerate}
\item
If $m\geq 6$ is even then
$$
\begin{array}{rcl}
N_2(\spin_{m}) = N_2^+ +N_2^- &=& \displaystyle
\sum_{\eta\in J_2(C)}h^0_+(\peta,m\xieta)
+ \sum_{\eta \not= 0}h^0_+(\peta^-,m\xieta), \\
N_2^+ -N_2^- &=& \displaystyle\sum_{\eta\in J_2(C)}h^0_-(\peta,m\xieta)
+ \sum_{\eta \not= 0}h^0_-(\peta^-,m\xieta).\\
\end{array}
$$
\item
If $m\geq 5$ is odd then
$$
\begin{array}{rcl}
N_2(\spin_{m}) = N_2^+ +N_2^- &=&\displaystyle \sum_{\eta\in
J_2(C)}h^0_+(\peta,m\xieta),\ \ \ \ \ \ \ \ \ \ \ \ \ \ \ \ \ \ \ \ \ \ \ \\
-N_2^+ +N_2^- &=& \displaystyle\sum_{\eta\in J_2(C)}h^0_-(\peta,m\xieta).\\
\end{array}
$$
\end{enumerate}
\end{theo}

\pf
Equations (\ref{spin6}), (\ref{spin2n}), (\ref{spin2n+1}) and (\ref{abtheta}).
\qed

\begin{rems}\rm
{\it (i)}
Of course, when $m$ is even $\peta$ and $\peta^-$ are interchangeable as far as
the dimensions alone are concerned. But in view of the constructions which
follow in the next sections, together with the remarks in section \ref{dynkin},
it will be seen that the right-hand side of part 1 is the correct way to write
these identities.

{\it (ii)}
For $m=3$ there is an analogous statement obtained by replacing $N_2$ by $N_4$.
This case has been discussed at length in \cite{OP}.

\end{rems}

These results imply, via (\ref{verlinde}), (\ref{di2}) and \ref{conj} (assuming
the validity of the latter) identities among $h^0$s which we can summarise in
the following table.

\begin{equation}
\label{dims}
\begin{array}{|c||l|l|}
\hline
&&\\
&h^0(\mm(\spin_m),\Theta(\c^m))&h^0(\mm^-(\spin_m),\Theta(\c^m))\\
&&\\
\hline
&&\\
{\rm odd}\ m& \displaystyle\sum_{\eta \in J_2} h^0_+(\peta, m \xieta)&
\displaystyle\sum_{\eta \in J_2} h^0_-(\peta, m \xieta)\\
&&\\
\hline
&&\\
{\rm even}\ m&\displaystyle\sum_{\eta \in J_2} h^0_+(\peta, m \xieta)&
\displaystyle\sum_{\eta \in J_2} h^0_-(\peta, m \xieta)\\
&&\\
 & \displaystyle +
\sum_{\eta \in J_2\backslash \{0\}} h^0_+(\peta^-, m \xieta)
& \displaystyle +
\sum_{\eta \in J_2\backslash \{0\}} h^0_-(\peta^-, m \xieta) \\
&&\\
\hline
\end{array}
\end{equation}

\bigskip

Although theorem \ref{numer} is stated only for $m\geq 5$, table (\ref{dims})
is in fact valid for {\it all} $m\in \n$, as we shall see next by a
case-by-case examination.

\begin{ex} $\bf m=1.$ \rm
We can identify $\mm(\spin_1)$ with $J_2(C)$; while $\mm^-(\spin_1)$ is empty.
So (\ref{dims}) is trivial in this case.
\end{ex}

\begin{ex} $\bf m=2.$ \rm
{}From example \ref{clifm=2}, $\mm^{\pm}(\spin_2)$ are both copies of
$\pic(C)$. However, the orthogonal theta bundle $\Theta(\c^2)$ is defined {\it
only} on the degree 0 component $J(C)$ of $\mm(\spin_2)$, where it can be
identified with the line bundle $8\theta$
where $\theta$ is a theta divisor on the Jacobian. So in this case again one
may readily check from (\ref{abtheta}) that the identity in (\ref{dims}) holds.
But in fact we have more:
$$
\begin{array}{rcl}
H^0_+(J(C),8\theta) &\cong& \displaystyle \sum_{\eta \in
J_2}H^0(\peta,2\xieta),\\
H^0_-(J(C),8\theta) &\cong& \displaystyle \sum_{\eta \in J_2\backslash
\{0\}}H^0(\peta^-,2\xieta).\\
\end{array}
$$

One can see this as follows.
Start with the fact that for any symmetric line bundle $\ll$ on an abelian
variety $A$ one has
$$
H^0(A,\ll^4) \cong \bigoplus_{\alpha\in \ahat}H^0(A,\ll\otimes \alpha)
$$
where $\ahat = \{\alpha \in \pic^0 A | \alpha^2 = \oo_A\}$. This isomorphism is
obtained by pulling back sections under the squaring map $[2]:A\rightarrow A$;
for each $\alpha$ we have $[2]^*(\ll\otimes \alpha) = \ll^4$. Noting that $[2]$
commutes with the involution $[-1] :A\leftrightarrow A$ we see that in fact
$
H^0_{\pm}(A,\ll^4) \cong \bigoplus_{\alpha\in \ahat}H^0_{\pm}(A,\ll\otimes
\alpha).
$

If $\theta$ represents a principal polarisation on $A$ then we can identify
$A_2 \cong \ahat$ by $\eta \mapsto \theta_{\eta} - \theta$ where $\theta_{\eta}
= t_{\eta}^*\theta$. In particular, when $A=J(C)$ and $\ll = 2\theta$ we obtain
$$
H^0_{\pm}(J,8\theta) \cong \bigoplus_{\eta\in J_2}H^0_{\pm}(J,\theta +
\theta_{\eta}).
$$
On the other hand $H^0_{\pm}(J,\theta + \theta_{\eta}) \cong
H^0(\peta^{\pm},2\xieta)$ for $\eta\not=0$. To see this choose $\zeta \in J(C)$
such that $\zeta^2 = \eta$. Then translation by $\zeta$ induces an isomorphism
$
t_{\zeta}^* : H^0(J,\theta + \theta_{\eta})\widetilde{\rightarrow}
H^0(J,2\theta),
$
which is equivariant with respect to $[-1]^*$ acting on the first space and
$t_{\eta}^*$ acting on the second. Hence the spaces
$H^0_{\pm}(J,\theta + \theta_{\eta})$ identify with the $\pm$-eigenspaces in
$H^0(J,2\theta)$ under $t_{\eta}^*$; and these in turn are classically
identified with $H^0(\peta^{\pm},2\xieta)$ by the Schottky-Jung relations.

\end{ex}

\begin{ex} $\bf m=3.$ \rm
We have seen (see example \ref{clifm=3}) that $\Theta(\c^3)$ restricts to
$\ll_0^4 \rightarrow \su(2)=\mm(\spin_3)$ and to $\ll_1^2 \rightarrow \su(2,1)
= \mm^-(\spin_3)$---in each case the anticanonical line bundle.
The coincidence of dimensions in (\ref{dims}) for this case was the first to be
observed, and has been explored in \cite{OP} and \cite{schottky}.
\end{ex}

\begin{ex}
\label{4theta}$\bf m=4.$ \rm
By example \ref{clifm=4} we have
$\mm(\spin_4) \cong \su(2)\times \su(2)$; while the 4-dimensional orthogonal
representation is the tensor product of the 2-dimensional representations of
the distinct factors, from which the theta bundle in (\ref{dims}) is
$\Theta(\c^4) = pr_+^*\ll_0^2 \otimes pr_-^* \ll_0^2$, where $\ll_0 = \Theta
(\c^2)$ is the ample generator of $\pic\ \su(2) = \z$. Thus
$$
H^0(\mm(\spin_4), \Theta(\c^4)) \cong \bigotimes^2 H^0(\su(2),\ll_0^2),
$$
and since $h^0(\su(2),\ll_0^2) = 2^{g-1}(2^g +1)$ we can at once verify the
case $m=4$ of the table. Moreover, just as for the case $m=2$ we can say more:
$$
\begin{array}{rcl}
S^2 H^0(\su(2),\ll_0^2) &\cong& \displaystyle \sum_{\eta \in
J_2}H^0_+(\peta,4\xieta),\\
\bigwedge^2 H^0(\su(2),\ll_0^2) &\cong& \displaystyle \sum_{\eta \in
J_2\backslash \{0\}}H^0_+(\peta^-,4\xieta).\\
\end{array}
$$

As in the case $m=2$ we have written down not just an identity of dimensions,
but in fact an isomorphism of vector spaces. In this case it follows easily
(for $C$ without vanishing theta-nulls) from the results of \cite{B2}.

Similarly $\mm^-(\spin_4)$ is isomorphic to $\su(2,1) \times \su(2,1)$, with
theta bundle $\Theta(\c^4) \cong pr_+^*\ll_1 \otimes pr_-^* \ll_1$, where
$\ll_1$ is the generator of $\pic\ \su(2,1)$. This time $h^0(\su(2,1),\ll_1^2)
= 2^{g-1}(2^g -1)$, the results of Beauville can be applied to give
isomorphisms
$$
\begin{array}{rcl}
S^2 H^0(\su(2,1),\ll_1) &\cong& \displaystyle \sum_{\eta \in
J_2}H^0_-(\peta,4\xieta),\\
\bigwedge^2 H^0(\su(2,1),\ll_1) &\cong& \displaystyle \sum_{\eta \in
J_2\backslash \{0\}}H^0_-(\peta^-,4\xieta),\\
\end{array}
$$
and again table (\ref{dims}) is verified in this case.

\end{ex}

\begin{ex}
\label{6theta}
$\bf m=6.$ \rm
In this case example \ref{clifm=6} identifies $\Theta (\c^6) \rightarrow
\mm^{\pm}(\spin_6)$ with $\ll_0^2 \rightarrow \su(4)$ and $\ll_2 \rightarrow
\su(4,2)$ respectively; and table (\ref{dims}) for the case $\mm^+(\spin_6)$
suggests isomorphisms
$$
\begin{array}{rcl}
H^0_+(\su(4),\ll_0^2) &\cong& \displaystyle \sum_{\eta \in
J_2}H^0_+(\peta,6\xieta),\\
H^0_-(\su(4),\ll_0^2) &\cong& \displaystyle \sum_{\eta \in J_2\backslash
\{0\}}H^0_+(\peta^-,6\xieta);\\
\end{array}
$$
for {\it some} $\pm$-decomposition of $H^0(\su(4),\ll_0^2)$. We shall observe
in remark \ref{dynkin} in the next section that the involution of
$\mm(\spin_6)$ analogous to the exchange of factors in example \ref{4theta} is
the dualising involution $E\mapsto E^{\vee}$ of $\su(4)$. We therefore
conjecture isomorphisms as above where $H^0_{\pm}$ on the left-hand side are
the $\pm$-eigenspaces for the dualising involution.

For $\mm^-(\spin_6) = \su(4,2)$, remark \ref{dynkin} will say that
the analogous involution is $E\mapsto E^{\vee}\otimes \nm E$,
where $\nm E$ is the fixed spinor norm, satisfying $(\nm E)^2 = \det E$.
Then we expect, in this case:
$$
\begin{array}{rcl}
H^0_+(\su(4,2),\ll_2) &\cong& \displaystyle \sum_{\eta \in
J_2}H^0_-(\peta,6\xieta),\\
H^0_-(\su(4,2),\ll_2) &\cong& \displaystyle \sum_{\eta \in J_2\backslash
\{0\}}H^0_-(\peta^-,6\xieta).\\
\end{array}
$$

\end{ex}

\section{Constructing the homomorphisms: the Jacobian}

Our aim in this section and the next is to construct homomorphisms:
\begin{equation}
\label{oddhomo}
H^0(\mm^{\pm}(\spin_{2n+1}), \Theta(\c^{2n+1}))^{\vee} \rightarrow
\sum_{\eta \in J_2}H^0_{\pm}(\peta,({2n+1})\xieta)
\end{equation}
\begin{equation}
\label{evenhomo}
\begin{array}{rcl}
H^0(\mm^{\pm}(\spin_{2n}), \Theta(\c^{2n}))^{\vee} &\rightarrow&
\displaystyle
\sum_{\eta \in J_2}H^0_{\pm}(\peta,2n\xieta) \\
&&\ \ \ \ \ \ \displaystyle
\oplus \sum_{\eta \not= 0}H^0_{\pm}(\peta^-,2n\xieta)\\
\end{array}
\end{equation}
Table (\ref{dims}) says that in each case the left- and right-hand sides have
the same dimension, and so we naturally conjecture that these homomorphisms are
isomorphisms, though we shall not prove this here.
As remarked in the last section, the first few cases $m=1,2,3,4$ are well
understood, and we do have natural isomorphisms (\ref{oddhomo}),
(\ref{evenhomo}) in these cases.

\subsection{The splitting for even spin groups}
\label{dynkin}

We shall at this point attempt to
explain the splitting $H^0(\Theta(\c^m)) = \sum_{\eta\in J_2} \oplus
\sum_{\eta\not=0} $ when $m$ is even. Note that for even $m$ the group
$\spin_m$ carries a unique nontrivial outer automorphism:
$$
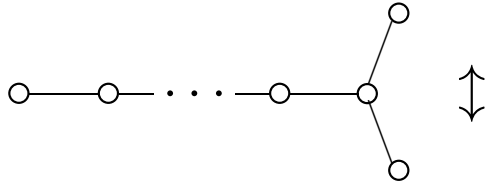
\dynk\dddnu......\updownarrow
$$
This induces an involution $\sigma : \mm(SC_m) \leftrightarrow \mm(SC_m)$ which
preserves the spinor norm and hence acts on $\mm(\spin_m)$ and
$\mm^-(\spin_m)$.
(More concretely, the group automorphism is obtained by Clifford conjugation by
a vector in $\c^m$; when $m$ is {\it odd} this is an inner automorphism.
$\sigma$ then acts on Clifford bundles by conjugating transition functions with
respect to some open cover; when $m$ is odd this action is still defined, of
course, but is trivial.)

On the other hand this automorphim preserves the orthogonal representation, and
hence the isomorphism class of the line bundle $\Theta(\c^m)$; accordingly
$H^0(\Theta(\c^m)$ splits into $\pm$-eigenspaces $H^0_+ \oplus H^0_-$
under the action of $\sigma$.

For even $m$ we now expect the following refinement of (\ref{evenhomo}):
\begin{equation}
\label{refine}
\begin{array}{rcl}
H^0_+(\mm^{\pm}(\spin_m), \Theta(\c^m))^{\vee} &\rightarrow & \displaystyle
\sum_{\eta \in J_2}H^0_{\pm}(\peta,m\xieta),\\
&&\\
H^0_-(\mm^{\pm}(\spin_m), \Theta(\c^m))^{\vee} &\rightarrow &\displaystyle
\sum_{\eta \not= 0}H^0_{\pm}(\peta^-,m\xieta).\\
\end{array}
\end{equation}

Note that this is exactly what happens in the case $m=4$ (see example
\ref{4theta}), where $\sigma$ simply exchanges factors in the products
$\mm(\spin_4) = \su(2) \times \su(2)$ and $\mm^-(\spin_4) = \su(2,1) \times
\su(2,1)$, so that the summands $H^0_{\pm}$ in (\ref{refine}) are precisely
$S^2 H^0$ and $\bigwedge^2 H^0$ appearing in \ref{4theta}.

In the case $m=6$---recall example \ref{clifm=6}---we have $\mm(\spin_6) =
\su(4)$ and $\mm^-(\spin_6) = \su(4,2)$. In each case not only the determinant
$\det E$ of vector bundles $E$ is fixed, but also
a square root $\nm E$ of $\det E$. Rank 4 vector bundles come from
$SC_6$-bundles via the first half-spinor representation $SC_6 \rightarrow
GL_4$; one easily checks that taking instead the second projection to $GL_4$
induces the involution on rank 4 vector bundles $\sigma : E\leftrightarrow
E^{\vee}\otimes \nm E$.

Thus in the case $m=6$ of (\ref{refine}) we expect the situation already
described in example \ref{6theta} of the previous section.

\subsection{The main construction}
\label{main}

We shall construct the homomorphisms (\ref{oddhomo}), (\ref{evenhomo}) one
summand at a time, concentrating in this section on the summand $\eta =0$, i.e.
the projection to Jacobian thetas (see (\ref{s0}) and corollary \ref{par0}).
In the next section section we shall construct the remaining projections to the
Prym thetas (corollary \ref{pareta}).

It will be convenient to denote the two-component variety $\mm(\spin_m) \cup
\mm^-(\spin_m)$ by $\nn(m)$. Let $\rho: SC_m \rightarrow SL(V)$ be a
representation satisfying condition (\ref{ss}) from section \ref{mod}.
(We are {\it mainly} concerned with the orthogonal representation $V=\c^m$, and
$r=m$
in proposition \ref{div} below.)

We consider the subset
\begin{equation}
\label{ddV}
\dd(V)  = \{(L,E)|H^0(C,L\otimes E(V)) \not= 0\} \subset J^{g-1}(C)\times
\nn(m).
\end{equation}
It is most important for us that in the case $V = \c^m$ is the vector
representation---by theorem \ref{gen0}---this is a proper subset in each
component of the product.

For the theta divisor $\theta$ in $J^{g-1}(C)$ and the theta bundle $\Theta(V)$
on $\nn(m)$ we shall abuse notation and use the same symbols to denote also
their pull-back to the product $J^{g-1}(C)\times \nn(m)$.

\begin{prop}
\label{div}
If $H^0(C,L\otimes E(V))$ is generically zero---as it is when $V=\c^m$---then
$\dd(V)$ is the support of a divisor in $|r\theta + \Theta(V)|$, where $r=\dim
V$.
\end{prop}

\pf
The representation $\rho$ induces
a morphism of varieties
$$
\alpha : \nn(m) \rightarrow \mm(SC_m) \rightarrow \mm(SL_r);
$$
and by definition $\dd(V)$ is (the support of) the pull-back via
$$
J^{g-1}(C) \times \nn(m) \map{id \times \alpha} J^{g-1}(C) \times \mm(SL_r)
\map{\otimes} \uu(r,r(g-1))
$$
(where, of course, we are identifying $\mm(SL_r)$ with the moduli space of
semistable vector bundles of rank $r$ and trivial determinant, and the second
map is tensor product of vector bundles)
of the canonical theta divisor (see section \ref{mod})
$$
\Theta_{r,r(g-1)} = \{F| H^0( F) \not= 0\} \subset \uu(r,r(g-1)).
$$
Since by hypothesis $H^0(C,L\otimes E(V))=0$ generically, it
follows that $\dd(V)$ is a well-defined divisor. The proposition now
follows from the discussion of section \ref{mod}: first of all, the pull-back
of $\Theta_{r,r(g-1)}$ to $J^{g-1}(C) \times \mm(SL_r)$ restricts on a fibre
$\{L\} \times \mm(SL_r)$ as the restriction of $\Theta_{r,0}$ from
$\uu(r,0)$---this is by definition of $\Theta_{r,0}$---and we have seen that
this is just $\Theta(\c^r)$. Hence the pull-back to $J^{g-1}(C) \times \nn(m)$
is $\Theta(V)$ on fibres $\{L\} \times \nn(m)$.

On the other hand, restriction to fibres $J^{g-1}(C) \times \{V\}$, for any
vector bundle $V$ with rank $r$ and trivial determinant, is well-known to be
independent of $V$---see for example \cite{OP} section 3.1.
Then it follows from (\ref{difr}) that $\Theta_{r,r(g-1)} $ restricts to
$r\theta$.
\qed

It follows from the K\"unneth theorem that $\dd(V)$ defines up to scalar a
tensor in
$$
H^0(\nn(m),\Theta(V))\otimes H^0(J^{g-1},r\theta),
$$
or equivalently a homomorphism, the projection of (\ref{oddhomo}),
(\ref{evenhomo}) at $\eta =0$,
\begin{equation}
\label{s0}
s_0:
H^0(\nn(m),\Theta(V))^{\vee} \rightarrow H^0(J^{g-1},r\theta).
\end{equation}
$s_0$ is dual to pull-back of hyperplane sections under the rational map
$
f_0 : \nn(m) \rightarrow |r\theta|
$
sending $E$ to $\dd(V)|_{J^{g-1}\times \{E\}}$.

Now suppose that the vector bundles $E(V)$ are self-dual. This happens when the
representation $V$ is symplectic or orthogonal---in particular if $V=\c^m$ or a
spin representation. Then it follows easily from Riemann-Roch and Serre duality
that the divisors $f_0(E)$ are {\it symmetric}:
$$
f_0 : \nn(m) \rightarrow |r\theta|_+ \cup |r\theta|_-
$$
where $|r\theta|_{\pm} = \bp H^0_{\pm}(J^{g-1},r\theta)$.

\begin{prop}
\label{f0}
When $V=\c^m$ is the standard orthogonal representation, for $m\geq 3$, $f_0$
respects parity:
$
f_0 : \mm^{\pm}(\spin_m) \rightarrow |m\theta|_{\pm}
$ respectively.
\end{prop}

Before proving this, we need to make a general remark about principally
polarised abelian varieties $(A,\Xi)$, where as usual $\Xi$ is a symmetric
divisor representing the polarisation. Let $\vartheta(A) = A_2 $ denote the set
of 2-torsion points, and let
$$
\begin{array}{rcl}
\vartheta^{+}(A) &=& \{ x\in A_2 \ |\  \mult_x \Xi \equiv 0\ \hbox{mod 2} \}\\
\vartheta^{-}(A) &=& \{ x\in A_2 \ |\  \mult_x \Xi \equiv 1\ \hbox{mod 2} \}\\
\end{array}
$$
In the case $(A,\Xi) = (J^{g-1} ,\theta)$ we shall write $\vartheta^{\pm}
(J^{g-1}) = \vartheta^{\pm}(C)$. These are the sets of even and odd theta
characteristics.

In the case of a Prym variety $(\peta ,\xieta)$ it is shown in \cite{OP},
proposition 2.3 that we can identify:
\begin{equation}
\label{Ptheta}
\begin{array}{rcr}
\vartheta^{\pm} (\peta) &=& \{\pi^* N = \pi^* (\eta \otimes N) \in \jtil  \\
&&\hbox{where}\ N, \eta\otimes N \in \vartheta^{\pm}(C)\}.\\
\end{array}
\end{equation}

On $\peta^-$ there is again an induced principal polarisation, though, as
remarked in section \ref{pryms},
no distinguished theta divisor. Thus we may talk about $\vartheta(\peta^-)$,
whose points are described (using the same methods as in \cite{OP}) by:
\begin{equation}
\label{P-theta}
\begin{array}{rcl}
\vartheta (\peta^-) &=& \{\pi^* N = \pi^* (\eta \otimes N) \in \jtil  \\
&&\hbox{where}\ N, \eta\otimes N \in \vartheta(C)\ \hbox{have {\it opposite}
parity}
\}.\\
\end{array}
\end{equation}
But for $\peta^-$ the partition into $\vartheta^{\pm}$ is no longer
well-defined.

\begin{lemm}
\label{basepts}
\begin{enumerate}
\item If $m$ is odd then $|m\Xi|_+$ (resp. $|m\Xi|_-$) is the linear subsystem
with
base-point set $\vartheta^-(A)$ (resp. $\vartheta^+(A)$).
\item If $m$ is even then $|m\Xi|_+$ is base-point free; while $|m\Xi|_-$ is
the linear subsystem with base-point set $A_2$.
\end{enumerate}
\end{lemm}

\pf See \cite{OP} lemma 2.2 (for part 1); or \cite{LB} chapter 4 section 7.
\qed

\medskip

{\it Proof of proposition \ref{f0}.}
We use proposition \ref{2.11}. Suppose first that $m$ is odd: then for $E\in
\nn(m)$ and any theta characteristic $L\in \vartheta(C)$
we have
$$
h^0(C,L\otimes E(\c^m)) \equiv h^0(C,L) + \deg\nm E
\quad {\rm mod}\ 2.
$$
So if $E\in \mm(\spin_m)$ then by definition $L\in f_0(E)$ for all {\it odd}
theta characteristics; so by the first part of the lemma $f_0(E) \in
|m\theta|_+$. Likewise $f_0(E) \in |m\theta|_-$ whenever $E\in \mm^-(\spin_m)$.

If $m$ is even then the same argument shows that $\mm^-(\spin_m)$ maps into
$|m\theta|_-$.
$\mm(\spin_m)$, on the other hand,
maps either into $|m\theta|_+$ or $|m\theta|_-$, and to see that it is not the
latter it suffices to exhibit a theta characteristic $\th$ and bundle $E\in
\mm(\spin_m)$ for which $H^0(C,\th\otimes E(\c^m)) =0$---which is possible by
theorem \ref{gen0th}.
\qed

\begin{cor}
\label{par0}
When $V=\c^m$ is the standard orthogonal representation the homomorphism $s_0$
respects parity:
$$
s_0 : H^0(\mm^{\pm}(\spin_m),\Theta(\c^m))^{\vee}
\rightarrow H^0_{\pm}(J^{g-1}, m\theta)
\qquad
{\sl respectively.}
$$
\end{cor}

\section{Constructing the homomorphisms: the Pryms}

We next explain how the construction of the previous section may be extended
to give the projections to the remaining summands, $\eta \not= 0$, in
(\ref{oddhomo}) and (\ref{evenhomo}). These are given in corollary
\ref{pareta}, which extends corollary \ref{par0}.

To begin, it is necessary to note that semistability of a bundle is preserved
under pull-back to the double covers.

\begin{lemm}
Let $p : \cctil \rightarrow C$ be any unramified cover of smooth projective
curves.  Then, if a vector bundle $V\rightarrow C$ is semistable then
$p^*V\rightarrow \cctil$ is semistable.
\end{lemm}

\pf By the Narasimhan-Seshadri theorem $V$ is induced from a projective unitary
representation of the fundamental group $\pi_1(C)$. Since $\cctil$ is an
unramified cover its fundamental group injects into $\pi_1(C)$, and the
restriction of the above representation then induces the pull-back bundle,
which is consequently semistable.
\qed

\begin{cor}
\label{5.5}
Let $p : \cctil \rightarrow C$ be as in the previous lemma, and $E\rightarrow
C$ a semistable $G$-bundle. Then $p^*E \rightarrow \cctil$ is semistable.
\end{cor}

\pf The same argument as in the above proof works for $G$-bundles by
Ramanathan's generalisation of the Narasimhan-Seshadri theorem \cite{R1};
alternatively apply the lemma to the adjoint bundle $\ad E$: by \cite{R2},
corollary 2.18, semistability of $E$ is equivalent to semistability of $\ad E$
as a vector bundle.
\qed

Let us now return to the double cover $p:\ctil \rightarrow C$. Noting that for
a Clifford bundle $E \rightarrow C$ the spinor norm satisfies $\nm ( p^*E) =
p^* \nm ( E)$, and this has even degree, it follows from corollary \ref{5.5}
that we obtain a morphism of moduli spaces
$$
u = p^*: \nn_C(m)  \rightarrow \mm_{\ctil}(\spin_m).
$$

\begin{prop}
\label{5.6} For any representation $SC_m \rightarrow SL(V)$ we have
 $u^*\Theta_{\mm_{\ctil}}(V) = 2\Theta_{\nn_C}(V)$.
\end{prop}

{\it Proof.}
Let $E\rightarrow C\times S$ be an arbitrary family of semistable
$SC_m$-bundles, and let $F=E(V)$ be the associated family of vector bundles via
the given representation. Let $\ftil = (p\times {\rm id})^*F$ be the pull-back
of the family by the double cover:
$$
\matrix{
\ctil \times S&&\map{p\times {\rm id}}&&C\times S\cr
&&&&\cr
&\hidewidth \pitil \searrow &&\swarrow \pi \hidewidth&\cr
&&&&\cr
&&S.&&\cr}
$$

It is clear from the discussion of section \ref{mod} that to prove the
proposition it suffices to show that
$$
\Theta(\ftil) = 2 \Theta (F):
$$
i.e. the line bundle $\Theta(V) \rightarrow \mm_C(SC_m)$ represents the functor
$E\mapsto \Theta (F)$, while the line bundle $u^*\Theta_{\mm_{\ctil}}(V)
\rightarrow \mm_C(SC_m)$ represents the functor $E\mapsto \Theta(\ftil)$.

So to compute $\Theta(\ftil)$, first note that by the projection formula
applied to $p\times {\rm id}$ we have, for $i=0,1$:
$$
\begin{array}{rcl}
R^i_{\pitil}(\ftil) &=& R^i_{\pi}(F\otimes p_* \oo_{\ctil}) \\
                    &=& R^i_{\pi}(F\oplus F\otimes \eta)\\
                    &=& R^i_{\pi}(F) \oplus R^i_{\pi}(F\otimes \eta).\\
\end{array}
$$
If we fit the direct images $R^i_{\pi}(F)$ into an exact sequence (\ref{1.3}),
and $R^i_{\pi}(F\otimes \eta)$ into a similar sequence with middle terms
${K^0}' \map{\phi'} {K^1}'$, then we get an exact sequence:
$$
0\rightarrow R^0_{\pitil}\ftil
\rightarrow K^0 \oplus {K^0}'
\map{
{\rm diag}(\phi,\phi')
} K^1 \oplus {K^1}'
\rightarrow R^1_{\pitil}\ftil\rightarrow 0.
$$
It follows at once that
$$
{\rm Det}(\ftil) = {\rm Det}(F)\otimes {\rm Det}(F\otimes \eta).
$$
But since the bundle $F$ has trivial determinant we can replace Det by $\Theta$
here. And since $\Theta(F\otimes \eta) = \Theta(F)$ by corollary \ref{1.6}, we
obtain $\Theta(\ftil) = 2 \Theta(F)$ as required.
\qed

We now consider
$$
\widetilde \dd(V) \subset \jtil \times
\mm_{\ctil}(\spin_m),
$$
denoting the divisor of (\ref{ddV}) and proposition \ref{div}
with $C$ replaced by $\ctil$.
As a consequence of proposition \ref{5.6}, we see that the pull-back via the
map
$$
 {\rm incl} \times u : (\peta \cup \peta^-) \times \nn(m) \rightarrow \jtil
\times
\mm_{\ctil}(\spin_m),
$$
of $\widetilde \dd(V)$ is---{\it provided $H^0(\ctil, L\otimes p^*E(V)) =0$ for
generic $(L,E) \in  \peta^{\pm}\times \mm^{\pm}(\spin_m)$}---a divisor
\begin{equation}
\label{eeV}
\ee_{\eta}(V) \in |2r \xieta + 2\Theta_{\nn}|.
\end{equation}

In the case of the orthogonal representation---except possibly for
$\mm^-(\spin_{2n})$---this is guaranteed by corollary \ref{serregen0}:

\begin{prop}
\label{gen0eta}
Let $V=\c^m$ be the standard orthogonal representation.
\begin{enumerate}
\item If $m$ is odd and $L\in \peta^-$ then $H^0(\ctil, L\otimes p^*E(V))
\not=0$ for all $E\in \mm^{\pm}(\spin_m)$.
\item In all other cases---except possibly $E\in \mm^-(\spin_{2n})$---we have
$H^0(\ctil, L\otimes p^*E(V)) =0$ for generic $(L,E) \in  \peta\times
\mm^{\pm}(\spin_m)$ or $\peta^-\times \mm^{\pm}(\spin_m)$.
\end{enumerate}
\end{prop}

\pf
Note that by corollary \ref{serre} (applied to $\ctil$), together with
(\ref{prym+-}), $h^0(\ctil, L\otimes p^*E(\c^m))$ is odd if $m$ is odd and
$L\in \peta^-$, hence nonzero as asserted. In all other cases, on the other
hand, $h^0(\ctil, L\otimes p^*E(\c^m))$ is even. So the proposition is
equivalent to \ref{serregen0}.
\qed

The next step is to observe that the divisor thus constructed has multiplicity
two:

\begin{prop}
Suppose that the representation $V$ is orthogonal, i.e. $SC_m \rightarrow
SO(V)$ for some invariant quadratic form on $V$. Then,
when
$\ee_{\eta}(V)$ is a divisor it
has multiplicity 2, i.e. $\ee_{\eta}(V) = 2 \dd_{\eta}(V)$ for a divisor
$\dd_{\eta}(V) \in |r \xieta + \Theta(V)|$.
\end{prop}

\pf We are here excluding the case $\nn(m) \times \peta^-$ for odd $m$. Then
we have already observed above that by corollary \ref{serre} and
(\ref{prym+-}), $h^0(\ctil, L\otimes p^*E(\c^m))$ is {\it even} for all $(L,E)
\in  \peta^{\pm}\times \mm^{\pm}(\spin_m)$ in the remaining cases.
The proposition then follows from the determinantal description of the functor
$\Theta$ given in section \ref{mod}.
\qed

Just as for the case $\eta =0$,
one now views the divisor $\dd_{\eta}(V)$ as giving a rational map---or, more
precisely, a {\it pair} of rational maps if $m$ is even:
$$
\begin{array}{rcl}
f_{\eta }^{\pm}: \nn(m) & \rightarrow& |m \xieta |_{\peta^{\pm}} \\
                E  &\mapsto & \{ L \in \peta^{\pm} | H^0(\ctil , L\otimes
p^*E(V)) \not= 0 \}.\\
\end{array}
$$
By Riemann-Roch, Serre duality and the fact that $E(V)$ is self-dual, one sees
again that $f_{\eta}^{\pm}(E)$ is a {\it symmetric} divisor.
This means that each component of $\nn(m)$ maps either into $|m\xieta
|_{+,\peta^{\pm}} =
\bp H^0_+(\peta^{\pm} , m \xieta)$ or into  $|m\xieta |_{-,\peta^{\pm}} = \bp
H^0_-(\peta^{\pm} , m \xieta)$; and the claim is:

\begin{prop}
\label{pmpm}
When $V=\c^m$ is the standard orthogonal representation, for $m\geq 3$, each
$f_{\eta} = f_{\eta}^+$ respects parity:
$$
f_{\eta}: \mm^{\pm}(\spin_m) \rightarrow |m\xieta |_{\pm,\peta};
$$
and when $m$ is even $f_{\eta}^-$ respects parity:
$$
f_{\eta}^-: \mm^{\pm}(\spin_m) \rightarrow |m\xieta |_{\pm,\peta^-}.
$$
\end{prop}

\pf
This uses proposition \ref{2.11} in the same way as the proof of proposition
\ref{f0}. We suppose first that $m$ is odd; choose $E\in \nn(m)$ and a theta
characteristic $L\in \vartheta(\peta)$. By (\ref{Ptheta}) this means $L= p^*N =
p^*(\eta\otimes N)$ for some theta characteristics $N$ and $\eta \otimes N$
{\it of the same parity}. Then we note that
$$
h^0(\ctil, L\otimes p^*E(\c^m)) = h^0(C,N\otimes E(\c^m)) + h^0(C,\eta \otimes
N\otimes E(\c^m)),
$$
and that (by \ref{2.11}) $h^0(C,N\otimes E(\c^m))$ is odd, and hence nonzero,
provided $E\in \mm(\spin_m)$ and $L\in \vartheta^-(\peta)$ {\it or} $E\in
\mm^-(\spin_m)$ and $L\in \vartheta^+(\peta)$. So by part 1 of lemma
\ref{basepts} it follows that $f_{\eta}$ maps $\mm^{\pm}(\spin_m)$ into
$|m\xieta|_{\pm,\peta}$ respectively.

If $m$ is even then the same argument, using part 2 of lemma \ref{basepts},
shows that $\mm^-(\spin_m)$ maps under $f_{\eta}$ into $|m\xieta|_{-,\peta}$;
and maps under $f_{\eta}^-$ into $|m\xieta|_{-,\peta^-}$ (see (\ref{P-theta})).

Finally consider $\mm(\spin_m) $ for even $m$: for any theta characteristics
$N, \eta \otimes N \in \vartheta(C)$, theorem \ref{gen0th} tells us that
$h^0(C,N\otimes E(\c^m)) = h^0(C,\eta \otimes N\otimes E(\c^m))=0$ for generic
$E\in \mm(\spin_m)$, so that generically $h^0(\ctil, L\otimes p^*E(\c^m)) = 0$
for any theta characteristic $L \in \vartheta(\peta)$ or $\vartheta(\peta^-)$.
This implies (for each of $\peta$, $\peta^-$) that $\mm(\spin_m) $ does not map
into $|m\xieta|_-$ and therefore maps into $|m\xieta|_+$ as required.
\qed

As a consequence, pull-back of hyperplane sections under each $f_{\eta}^{\pm}$
is dual to a homomorphism (the analogue of (\ref{s0}))
\begin{equation}
\label{seta}
s_{\eta}^{\pm}:
H^0(\nn(m),\Theta(V))^{\vee} \rightarrow H^0(\peta^{\pm},r\xieta),
\end{equation}
where $r =\dim V$, and where $s_{\eta}^-$ is defined only for even $m$.
Moreover, proposition \ref{pmpm} says that
these homomorphisms {\it respect parity}---in other words
the analogue of corollary \ref{par0} for this situation is:

\begin{cor}
\label{pareta}
When $V=\c^m$ is the standard orthogonal representation the homomorphisms
$s_{\eta}^{\pm}$ respect parity:
$$
s_{\eta} = s_{\eta}^+ : H^0(\mm^{\pm}(\spin_m),\Theta(\c^m))^{\vee}
\rightarrow H^0_{\pm}(\peta, m\xieta)
\qquad
{\sl respectively,}
$$
and if $m$ is even then additionally:
$$
s_{\eta}^- : H^0(\mm^{\pm}(\spin_m),\Theta(\c^m))^{\vee}
\rightarrow H^0_{\pm}(\peta^-, m\xieta)
\qquad
{\sl respectively.}
$$
\end{cor}

\bigskip

{\addressit

 Department of Mathematical Sciences,

 Science Laboratories,

 South Road,

 Durham DH1 3LE,

 U.K.

}

\medskip

 {\addressit E-mail:} {\eightrm w.m.oxbury@durham.ac.uk  }

\end{document}